# Sur les zéros des fonctions $L$ automorphes de grand niveau

E. Kowalski, Rutgers University,    P. Michel, Université d'Orsay

July 31, 1997

**Résumé**

We study, on average over $f$, zeros of the $L$-functions of primitive weight two forms of level $q$ (fixed). We prove, on the one hand, density theorems for the zeros (similar to the results of Bombieri, Jutila, Motohashi, Selberg in the case of characters), which are applied in [KM1] to obtain a sharp unconditionnal estimate of the (analytic) rank of the new part of $J_0(q)$; and, on the other hand, non-vanishing theorems at the critical point, showing that a positive proportion of $L$-functions are non zero there.

## Table des matières

## 1 Introduction

Dans l'article compagnon [KM1] nous étudions, par la méthode des formules explicites de Mestre, le rang analytique de la partie nouvelle $J_0(q)^+$ de la jacobienne de la courbe modulaire $X_0(q)$. Pour s'affranchir de l'hypothèse de Riemann, suivant l'approche utilisée par Perelli et Pomykala [P-P] pour les tordues quadratiques des courbes elliptiques, nous sommes amenés à démontrer des théorèmes de densité pour les zéros des fonctions $L$ des formes automorphes de poids 2 et de niveau $q$.

De tels résultats ont été obtenus sous diverses formes pour les caractères de Dirichlet. Nous établissons ici trois résultats dont l'efficacité varie suivant les circonstances. Deux d'entre eux sont utilisés dans [KM1] mais tous pourront peut-être trouver d'autres applications.

Rappelons nos notations. Pour toute famille de formes automorphes $(f_i)_i$ et tout complexes $(\alpha_f)_f$, nous noterons

$$\sum_f^h \alpha_f := \sum_f \frac{1}{4\pi(f,f)} \alpha_f$$

($h$ pour "harmonique".) Ici $(f,f)$ est le produit scalaire de Petersson, normalisé de la manière suivante

$$(f,g) = \int_{\Gamma_0(q)\backslash \mathbf{H}} f(z)\overline{g(z)} dx dy$$

pour $f \in S_2(q)$, $g \in S_2(q)$ des formes paraboliques de poids 2.

Nous écrivons le développement de Fourier d'une forme parabolique sous la forme

$$f(z) = \sum_{n\geq 1} \lambda_f(n) n^{1/2} e(nz)$$

de sorte que la fonction $L$ de $f$ est

$$L(f,s) = \sum_{n\geq 1} \lambda_f(n) n^{-s}$$

et que sa droite critique est alors $\operatorname{Re}(s) = 1/2$.

Pour toute forme $f$, nous notons $N(f, \alpha, T, T')$ le nombre de zéros $\beta + i\gamma$ de $L(f)$ qui vérifient

$$\beta \geq \alpha$$
$$T \leq \gamma \leq T'$$

et $N(f, \alpha, T) = N(f, \alpha, -T, T)$.

En premier lieu, nous prouvons l'analogue du théorème 18 de Bombieri [Bo].

**Théorème 1.1** *Il existe $A > 0$ et $B > 0$ tels que pour tout $\alpha$, $\frac{1}{2} \leq \alpha \leq 1$ et $T \geq 1$ on ait*

$$\sum_{f \in S_2(q)^+}^h N(f, \alpha, T) \ll \frac{1}{q} T^A q^{\frac{3(1-\alpha)}{2-\alpha}} (\log q)^B$$

*d'autre part si $q$ est premier on a*

$$\sum_{f \in S_2(q)^+} N(f, \alpha, T) \ll T^A q^{\frac{3(1-\alpha)}{2-\alpha}} (\log q)^{B+1}.$$

La deuxième inégalité resulte aussitôt de la première grâce à la minoration $(f,f) \gg q/\log q$ valable pour $q$ premier [HL, GHL]. Dans le cas général, on a le corollaire suivant sans poids $1/4\pi(f,f)$ qui est plus faible dans l'exposant en $q$ mais qui est suffisant pour une des applications de [KM1]:

**Corollaire 1.1** *Il existe $A > 0$ et $B > 0$ et $c > 0$ tels que pour tout $\alpha$, $\frac{1}{2} \leq \alpha \leq 1$ et $T \geq 1$ on ait*

$$\sum_{f \in S_2(q)^+} N(f, \alpha, T) \ll T^A \dim J_0^n(q) q^{-c(\alpha - 1/2)} (\log q)^B.$$

Au voisinage de l'abcisse de convergence absolue, le facteur logarithmique ci-dessus devient prédominant. Le résultat suivant le supprime,

**Théorème 1.2** *Il existe une constante positive $B$ telle que, pour $1 - 1/13 \leq \alpha \leq 1$, $T \geq 1$ et tout $\varepsilon > 0$, on a*

$$\sum_{f \in S_2(q)^+} N(f, \alpha, T) \ll T^B q^{(13+\varepsilon)\frac{1-\alpha}{2\alpha - 1}}.$$

Dans le cas des caractères de Dirichlet, le premier théorème de ce type est dû à Fogels [Fo], encore qu'il faille citer les travaux fondamentaux de Linnik sur le plus petit nombre premier dans une progression arithmétique [Li]. La contribution décisive suivante est due à Selberg qui, en introduisant la méthode des pseudo-caractères, à non seulement considérablement simplifié les preuves mais aussi amélioré les estimations de densité pour leur donner la forme moderne suivante [M-S] (avec des notations évidentes que nous ne rappelons pas)

$$\sum_{\chi \,(mod\ q)} N_\chi(\alpha, T) \ll_\epsilon (qT)^{(3+\epsilon)(1-\alpha)}.$$

Cette majoration a été améliorée par Motohashi [Mot] puis Jutila [Ju], qui détient le record actuel avec la constante 3 remplacée par 2 dans l'expression précédente (et avec la condition $\alpha \geq 4/5$). Il est probable que la méthode présentée ici peut servir de modèle pour montrer des Théorèmes de densité analogues pour les zéros de fonctions $L$ de puissance symétriques de formes modulaire ou plus généralement de représentations automorphes sur $GL(n)$.

Pour intéressants qu'ils soient ces Théorèmes sont insuffisants pour les applications arithmétiques que nous avons en vue. Nous nous restreignons dans la suite au cas où $q$ est premier: le résultat suivant compte les zéros dans de petits intervalles verticaux avec une puissance de $\log q$ optimale pour la technologie actuelle; il est particulièrement efficace au voisinage de la droite critique (voir l'application dans [KM1] où on en déduit la borne optimale, à la constante près, $\text{rang}_a J_0(q) = O(\dim J_0(q))$; ce Théorème de densité correspond au théorème 4 de Selberg [Sel]:

**Théorème 1.3** *Soit $q$ premier. Il existe des constantes absolues $A > 0$ et $\kappa > 0$ telles que pour tout $T \geq 1$, et tout couple $(t_1, t_2)$ de réels satisfaisant*

$$-T \leq t_1 < t_2 \leq T$$

$$t_2 - t_1 \geq \frac{\kappa}{\log q}$$

*et pour tout $\alpha \geq \frac{1}{2} + \frac{1}{\log q}$, et $0 < c < 1/8$ on a*

(1) $$\sum_{f \in S_2(q)^+}^h N(f, \alpha, t_1, t_2) \ll T^A q^{-c(\alpha - \frac{1}{2})} (\log q)(t_2 - t_1)$$

*et de même*

(2) $$\sum_{f \in S_2(q)^+} N(f, \alpha, t_1, t_2) \ll T^A q^{1-c(\alpha-\frac{1}{2})} (\log q)(t_2 - t_1).$$

La preuve de ces deux majorations est similaire à la technique employée par Selberg; en particulier, on doit donner de bonnes majorations pour la moyenne quadratique

$$\sum_{f \in S_2(q)^+}^h |L(f,s)M(f,s)|^2,$$

où $M(f,s)$ est un certain "mollifieur" et $s$ est proche de la droite critique. Comme nous l'a fait remarqué H. Iwaniec de telles majorations peuvent être utilisées pour montrer des théorèmes de non-annulation pour les fonctions $L(f,s)$ au point critique $s = 1/2$. Suivant cette suggestion pour montrons le Théorème de non-annulation suivant qui est l'analogue "suivant le niveau" de Théorèmes de non-annulation "suivant le poids" de Iwaniec-Sarnak [IS]:

**Théorème 1.4** *Pour tout $\epsilon > 0$, et pour $q$ premier assez grand, on a les minorations*

(3) $$\sum_{\substack{f \in S_2(q)^+ \\ L(f,1/2) \neq 0}}^h 1 \geq \frac{1}{6} - \epsilon.$$

(4) $$|\{f \in S_2(q)^+, \ L(f,1/2) \neq 0\}| \geq (\frac{1}{6} - \epsilon)|S_2(q)^+|$$

Remarquons que si $q$ est premier l'application $f \to 1/4\pi(f,f)$ définit asymptotiquement une mesure de probabilité sur $S_2(q)^+$ on a en effet

$$\sum_f^h 1 = 1 + O(q^{-3/2});$$

les minorations (3) et (4) signifient donc que pour la mesure ainsi définie ou pour la mesure "naturelle" $f \to 1/|S_2(q)^+|$, il y a une proportion positive de formes $f$ telles que $L(f, 1/2) \neq 0$. Remarquons d'autre part que les formes primitives $f$ associées à la valeur propre -1 de l'involution d'Atkin-Lehner représentent approximativement 50 % de l'ensemble $S_2(q)^+$ (pour l'une ou l'autre des deux mesures) et que pour celles-ci on a $L(f, 1/2) = 0$, la minoration (3) (resp. (4)) signifie donc qu'au moins 33% des fonctions $L$ des formes primitives associées à la valeur propre 1 de l'involution d'Atkin-Lehner ne s'annulent pas au point critique relativement à la mesure harmonique ou à la mesure naturelle (la proportion de 1/3 est en parfaite analogie avec un résultat de non-annulation pour les caractères de Dirichlet démonté par Iwaniec – non publié–).

La minoration (4) peut encore se réécrire, suivant la terminologie de Merel [Me], en terme du "quotient d'enroulement":

$$\dim J_0^e(q) \geq (\frac{1}{6} - \epsilon) \dim J_0(q).$$

Suivant Duke [Du2], on a également le corollaire suivant:

**Corollaire 1.2** *Soit $M_C^*$ l'espace de Kohnen des formes modulaires de poids $3/2$ pour $\Gamma_0(4q)$ (ie. les formes dont le n-ième coefficient de Fourier s'annule sauf si $-n \equiv 0, 1 \pmod 4$ et $(\frac{-n}{q}) \neq 1$) et $\Theta_q$ le sous-espace engendré par les fonctions Théta; on a alors pour tout $\epsilon > 0$ et $q$ premier assez grand la minoration*

$$\dim \Theta_q \geq (\frac{1}{3} - \epsilon) \dim M_C^*.$$

( Duke avait obtenu la minoration plus faible $\dim \Theta_q \gg \dim M_C^* / \log^2 q$).

Disons enfin que par des méthodes similaires, on peut aussi considérer les formes primitives $f$ associées à la valeur propre -1 de l'involution d'Atkin-Lehner (on a alors $L(f, 1/2) = 0$ ) et de montrer que pour une proportion positive de celles-ci, la dérivée $L'(f, 1/2)$ au point critique est non-nulle, nous reviendrons sur ce point dans [KM2].

**Remarque**. — Comme nous l'a fait remarquer H. Iwaniec, il est possible de remplacer la constante $1/6$ par $1/4$ dans les minorations (3) et (4), pour renvoyons pour cela à [IS]

**Remarque**. — La restriction à $q$ premier dans ce dernier énoncé, intervient du fait de l'utilisation de la formule de Petersson qui ne s'applique directement qu'à une base orthonormée complète de l'espace des formes automorphes et non seulement à la famille des formes primitives de Hecke. Pour $q$ premier, cependant, du fait de l'absence de formes paraboliques de poids 2 de niveau 1 (car $X_0(1)$ est de genre 0), toutes les formes de poids 2 sont nouvelles. Il est probable que par un argument d'inclusion-exclusion sur les formes anciennes on puisse traiter de même le cas général

Le plan de cet article est le suivant, après une section consacrée à diverses généralités nous montrons dans les sections 3 et 4 les inégalités (1) et (3), dans les deux suivantes 5 et 6 nous montrons (2) et (4) (cette organisation nous a paru plus naturelle dans la mesure où les résultats en moyenne harmonique découlent plus directement de la formule de Petersson alors que les résultats en moyenne "naturelle" sont une modification un peu technique des inégalités (1) et (3)). Dans la section 7 nous montrons le Théorème 1.1 et nous terminons par la preuve du Théorème 1.2.

**Remerciements:** Cet article doit beaucoup à Henryk Iwaniec et il serait vain d'énumérer toutes les améliorations qu'il a apporté tant du point de vue de la clarification des preuves que des mathématiques, disons qu'en particulier que la section 6 est largement inspirée par ses idées. Nous le remercions également ainsi que Peter Sarnak pour nous avoir montré leur travail [IS]. Nous remercions encore Bill Duke et Etienne Fouvry pour de stimulantes discussions et Matti Jutila pour les références qu'il nous a donné.

## 2 Préliminaires

On note $S_2(q)$ (resp. $S_2(q)^+$) l'espace des formes modulaires paraboliques de poids 2 et de niveau $q$ (resp. l'ensemble des formes primitives.)

Toute forme parabolique $f$ a un développement de Fourier (à la pointe infinie)

$$f(z) = \sum_{n \geq 1} \lambda_f(n) n^{1/2} e(nz)$$

(comme d'habitude, $e(z) = \exp(2i\pi z)$.)

Soit $f \in S_2(q)^+$ une forme primitive, elle est alors normalisée par $\lambda_f(1) = 1$. Ses coefficients de Fouriers $\lambda_f(n)$ sont réels et vérifient les relations de Hecke

$$\lambda_f(l_1)\lambda_f(l_2) = \sum_{d|(l_1,l_2)} \varepsilon(d)\lambda_f\Big(\frac{l_1 l_2}{d^2}\Big).$$

La fonction $L$ associée à $f$ est

$$L(f,s) = \sum_{n \geq 1} \lambda_f(n) n^{-s} = \prod_p (1 - \lambda_f(p) p^{-s} + \varepsilon(p) p^{-2s})^{-1}$$

où $\varepsilon$ est le caractère de Dirichlet trivial modulo $q$.

Cette fonction admet un prolongement analytique comme fonction entière sur tout le plan complexe et la fonction $L$ complétée

$$\Lambda(f,s) = \Big(\frac{q}{2\pi}\Big)^{\frac{s}{2}} \Gamma(s+\frac{1}{2}) L(f,s)$$

a l'équation fonctionnelle

$$\Lambda(f,s) = \varepsilon_f \Lambda(f, 1-s)$$

avec $\varepsilon_f = q^{1/2} \lambda_f(q) = \pm 1$ [Miy].

## 2.1 Base orthonormée, sommes de Kloosterman et grand crible

Soit $\mathcal{F}$ est une base orthonormée de $S_2(q)$, la "formule des traces" de Petersson relie les coefficients de Fourier des éléments de $\mathcal{F}$ à des sommes de sommes de Kloosterman (cf. [DFI]): pour tout $m, n \geq 1$ on a l'égalité

(5) $$\frac{1}{4\pi} \sum_{f \in \mathcal{F}} \lambda_f(m)\lambda_f(n) = \delta_{m,n} - 2\pi \sum_{c \geq 1} \frac{S(m,n;cq)}{cq} J_1\Big(\frac{4\pi\sqrt{mn}}{cq}\Big)$$

où $\delta_{m,n}$ est le symbole de Kronecker,

$$S(m,n;c) = \sum_{\substack{x \bmod c \\ (x,c)=1}} e\Big(\frac{mx + n\overline{x}}{c}\Big)$$

est la somme de Kloosterman et $J_1(x)$ est la fonction de Bessel d'ordre 1. Nous utiliserons les majorations suivantes de $S(m,n;c)$ et $J_1(x)$:

$$S(m,n;c) \leq (m,n,c)^{1/2} c^{1/2} \tau(c)$$

$$J_1(x) \ll \min(x, x^{-1/2}).$$

On en déduit pour tout $m$ et $n$

(6) $$\frac{1}{4\pi} \sum_{f \in \mathcal{F}} \lambda_f(m)\lambda_f(n) = \delta_{m,n} + O((m,n,q)^{1/2}(mn)^{1/2} q^{-3/2}).$$

Duke, Friedlander et Iwaniec ont combiné cette formule de trace avec le grand crible arithmétique pour donner des inégalités de grand crible quasi-optimales sur les coefficients de Fourier des formes $f \in \mathcal{F}$ ([DFI] Theorem 1.):

**Proposition 2.1** *Soit $\beta = (\beta_n)$ une suite de nombres complexes. Alors on a la majoration*

$$\sum_{f \in \mathcal{F}} \Big| \sum_{n \leq N} \beta_n \lambda_f(n) \Big|^2 \ll \Big(1 + \frac{N}{q}\Big) \sum_{n \leq N} |\beta_n|^2$$

*la constante sous-entendue étant absolue.*

**Remarque.** — Un facteur $\log N$ est présent dans [DFI], mais il peut être enlevé en utilisant une méthode due à Iwaniec (voir [Du1] Prop. 1 pour le cas des formes de poids 1.)

## 2.2  L'inverse et le carré de $L(f)$

Soit $f$ primitive de niveau $q$. Le produit eulérien pour $L(f,s)$ est

$$L(f,s) = \prod_p (1 - \lambda_f(p)p^{-s} + \varepsilon(p)p^{-2s})^{-1}.$$

**Lemme 2.1** $L(f)^{-1}$ *est donnée par la série de Dirichlet*

$$\begin{aligned}L(f,s)^{-1} &= \prod_p (1 - \lambda_f(p)p^{-s} + \varepsilon(p)p^{-2s}) \\ &= \sum_{l,m \geq 1} \varepsilon(m)\mu(l)|\mu(lm)|\lambda_f(l)(lm^2)^{-s}.\end{aligned}$$

*Preuve.* — La forme de la série s'obtient par multiplicativité, en notant que tout $n \geq 1$ s'écrit de façon unique $n = lm^2r$, avec $l$, $m$, $r$ premiers deux à deux, $l$ et $m$ sans facteurs carrés et chaque facteur premier de $r$ de valuation au moins 3.

□

**Lemme 2.2** $L(f)^2$ *est donnée par la série de Dirichlet*

$$L(f,s)^2 = \sum_{l,m \geq 1} \varepsilon(m)\tau(l)\lambda_f(l)(lm^2)^{-s}.$$

*Preuve.* — D'après les relations de multiplicativité, le coefficient de $n^{-s}$ dans $L(f,s)^2$ est

$$\begin{aligned}\sum_{d|n} \lambda_f(d)\lambda_f(n/d) &= \sum_{d|n} \sum_{m|(d,n/d)} \varepsilon(m)\lambda_f(n/m^2) \\ &= \sum_{m^2|n} \varepsilon(m)\lambda_f(nm^{-2}) \sum_{l|nm^{-2}} 1\end{aligned}$$

ce qui donne le résultat en sommant sur $n$.

□

La série donnant $L(f)^{-1}$ est absolument convergente là où celle de $L(f)$ l'est, de même pour $L(f)^2$.

## 3 Théorème de densité en moyenne harmonique

Celui-ci étant, du point de vue de l'application de [KM1], le plus intéressant, et également le plus difficile de toute façon, nous commençons par montrer le théorème 1.3, ou plus précisément la majoration (1) suivant l'approche de Selberg dans [Sel]. Celle-ci est basée sur l'étude de la somme

$$S = \sum_{f \in S_2(q)^+}^{h} |L(f, \beta + it)M(f, \beta + it)|^2$$

où $\beta = 1/2 + \delta$ et $\delta$ est réel, dans le cas le plus défavorable (celui, donc, où ce dernier théorème de densité est plus efficace que les précédents) est de l'ordre de $\frac{1}{\log q}$.

Cependant il est parfaitement possible de mener l'analyse pour toute valeur de $\delta$, $0 \leq \delta < 1/2$ mais nous l'éviterons par l'utilisation d'une astuce de convexité. Le cas $\delta = 0$ (correspondant aux zéros sur la droite critique) est intéressant pour le section (3) et nous le considérerons également, même lorsque son traitement divergera du cas $\delta \gg \frac{1}{\log q}$.

*Dans cette section et les section 4, 5 et 6, nous supposons que $q$ est un nombre premier.*

### 3.1 La somme principale

Pour $f \in S_2(q)^+$, on définit la fonction entière

$$Z(f, s) = \Lambda(f, s+it)\Lambda(f, s-it),$$

et puisque les coefficients de Fourier de $f$ sont réels nous aurons

(7) $$Z(f, \beta + it) = |\Lambda(f, \beta + it)|^2$$

et $Z(f)$ vérifie l'équation fonctionnelle simple

$$Z(f, s) = Z(f, 1-s).$$

### 3.2 Lemmes techniques

Nous fixons $N \geq 1$ un entier grand mais fixé, et choisissons (arbitrairement) un polynôme *à coefficients réels* $G(s)$ qui satisfasse:

(8) $$G(1-s) = G(s)$$

(9) $$G(-1/2) = G(-3/2) = \ldots = G(-N+1/2) = 0, G(1/2) = 1$$

On pose de plus

(10) $$H(s) = \Gamma(s + \frac{1}{2} + it)\Gamma(s + \frac{1}{2} - it)G(s+it)G(s-it).$$

et on construit successivement les fonctions suivantes:

(11) $$W_\delta(y) = \frac{1}{2i\pi} \int_{(2)} H(s+\beta)y^{-s}s^{-1}ds$$

$$
(12) \qquad V(y) = y^{-\delta}W_\delta(y) + y^\delta W_{-\delta}(y)
$$

$$
(13) \qquad U(y) = \sum_{d \geq 1} \frac{\varepsilon(d)}{d} V(yd^2).
$$

Le lemme suivant rassemble les bornes et majorations requises dans la suite.

**Lemme 3.1** *Les fonctions $W_\delta$ et $V$ vérifient, pour certaines constantes $B > 0$ dont la valeur peut varier d'une ligne à l'autre*

$$
(14) \qquad W_\delta(y),\ y^{-\delta}V(y),\ U(y) \ll y^{-1}(|t|+1)^B e^{-\pi|t|},\ y \to +\infty
$$

$$
(15) \qquad W_\delta(y) = H(\beta) + O(y^1(|t|+1)^B e^{-\pi|t|}),\ y \to 0
$$

*et, si $\delta \neq 0$, en écrivant $\zeta_q(s) = L(s, \varepsilon)$, au voisinage de $y = 0$*

$$
(16)\ V(y) = H(\beta)\zeta_q(1+2\delta)y^{-\delta} + H(1-\beta)\zeta_q(1-2\delta)y^\delta + O(y^{1/2}(|t|+1)^B e^{-\pi|t|}),\ y \to 0
$$

*alors que si $\delta = 0$,*

$$
(17) \qquad U(y) = -\frac{\phi(q)}{q}\log y + C_q + O(y^{1/2}(|t|+1)^B e^{-\pi|t|}),\ y \to 0
$$

*où*

$$
C_q = c_0 + \frac{c_1 + c_2 \log q}{q}
$$

*et $c_0, c_1, c_2$ sont des constantes ne dépendant que de $G$.*

Notons également, puisqu'on aura l'occasion de diviser par $H(\beta)$, que la formule de Stirling entraîne

$$
H(\beta) \gg e^{-\pi|t|}.
$$

*Preuve.* — Tout cela est très simple par transformation de Mellin et déformation de contour, en utilisant au besoin (9).

Par exemple, pour (16), par intégration complexe et (13), (12), (11) on a

$$
\begin{aligned}
U(y) &= \frac{1}{2i\pi}\int_{(2)} H(s+\beta)\zeta_q(1+2s+2\delta)y^{-s-\delta}s^{-1}ds \\
&\quad + \frac{1}{2i\pi}\int_{(2)} H(s+1-\beta)\zeta_q(1+2s-2\delta)y^{-s+\delta}s^{-1}ds
\end{aligned}
$$

et l'on déforme le contour jusqu'à la droite de partie réelle $-1/2$, ce qui amène à passer trois pôles en $s = \delta, 0, -\delta$, qui sont bien distincts par hypothèse. En $\delta$ et $-\delta$, les résidus sont respectivement

$$
\frac{1}{\delta}H(1/2)(1-q^{-1})
$$

et

$$
-\frac{1}{\delta}H(1/2)(1-q^{-1})
$$

dont la somme est nulle. Estimant l'intégrale obtenue sur $\mathrm{Re}(s) = -1/2$, le résidu en $s = 0$ fournit le lemme. Si $\delta = 0$,

$$U(y) = \frac{2}{2i\pi} \int_{(2)} H(s+1/2)\zeta_q(1+2s)y^{-s}s^{-1}ds,$$

si bien que cette fois-ci l'intégrande à un pôle double en 0 dont on calcule le résidu.

Les autres égalités sont similaires.

□

On rencontrera aussi la fonction arithmétique $\eta_t$

(18) $$\eta_t(n) = \sum_{l_1 l_2 = n} (l_1/l_2)^{it}.$$

**Lemme 3.2** *La fonction arithmétique $\eta_t$ est à valeurs réelles. De plus $\eta_t(n)$ est la valeur propre de l'opérateur de Hecke $T(n)$ pour une série d'Eisenstein (non-holomorphe) de niveau 1. En particulier, on a les relations*

(19) $$\eta_t(n)\eta_t(m) = \sum_{d|(n,m)} \eta_t\left(\frac{nm}{d^2}\right)$$

(20) $$\eta_t(nm) = \sum_{d|(n,m)} \mu(d)\eta_t\left(\frac{n}{d}\right)\eta_t\left(\frac{m}{d}\right)$$

*ainsi que la majoration*
(21) $$\eta_t(n) \ll \tau(n).$$

*Enfin on a les factorisations suivantes des séries de Dirichlet*

(22) $$\sum_{n \geq 1} \frac{\eta_t(n)}{n^s} = \zeta(s-it)\zeta(s+it),$$

(23) $$\sum_{n \geq 1} \frac{\eta_t(n)^2}{n^s} = \frac{\zeta(s-2it)\zeta(s)^2\zeta(s+2it)}{\zeta(2s)},$$

(24) $$\sum_{n \geq 1} \frac{\eta_t(n^2)}{n^s} = \frac{\zeta(s-2it)\zeta(s)\zeta(s+2it)}{\zeta(2s)},$$

*Preuve.* — Que $\eta_t(n)$ soit valeur propre de l'opérateur de Hecke est bien connu. La première relation est la relation de Hecke comme ci-dessus, et la seconde est son inversion par Möbius. Les factorisations des séries de Dirichlet, correspondent à celle des fonctions $L$ respectivement, de la forme elle-même, de sa convolution de Rankin-Selberg et de son carré symétrique. Elles se démontrent de manière élémentaire.

□

### 3.3 Expression de $|L(f, \beta + it)|^2$

On considère maintenant

$$
\begin{aligned}
I_\beta(f) &= \frac{1}{2i\pi} \int_{(2)} Z(f,s) G(s+it) G(s-it) \frac{ds}{s-\beta} \\
&= \frac{1}{2i\pi} \int_{(2)} L(f, s+it) L(f, s-it) H(s) \Big(\frac{q}{2\pi}\Big)^{2s} \frac{ds}{s-\beta}
\end{aligned}
$$

(voir (10)).

Puisque l'intégrande est à décroissance rapide sur les bandes verticales (grâce aux fonctions $\Gamma$), on peut déplacer le contour d'intégration jusqu'à la droite de partie réelle $-1/2$. Chemin faisant un pôle simple est obtenu en $s = \beta$. Appliquant à l'intégrale ainsi obtenue le changement de variable $s \mapsto 1-s$ et l'équation fonctionnelle de $Z(f)$ ainsi que (8), on voit que celle-ci n'est autre que $-I_{1-\beta}(f)$. Calculant le résidu, il vient

(25) $$\Big(\frac{q}{2\pi}\Big)^\beta H(\beta) L(f, \beta+it) L(f, \beta-it) = I_\beta(f) + I_{1-\beta}(f).$$

Mais on peut également exprimer $I_\beta(f)$ en développant $L(f)$ en série de Dirichlet

$$I_\beta(f) = \Big(\frac{q}{2\pi}\Big)^\beta \sum_{l_1, l_2 \geq 1} \lambda_f(l_1) \lambda_f(l_2) (l_1 l_2)^{-\beta} (l_1/l_2)^{it} W_\delta\Big(\frac{2\pi l_1 l_2}{q}\Big)$$

(voir (11)). Insérant cela dans (25) on obtient (voir (12))

$$\Big(\frac{q}{2\pi}\Big)^\delta H(\beta) L(f, \beta+it) L(f, \beta-it) = \sum_{l_1, l_2 \geq 1} \lambda_f(l_1) \lambda_f(l_2) (l_1 l_2)^{-1/2} (l_1/l_2)^{it} V\Big(\frac{2\pi l_1 l_2}{q}\Big)$$

que l'on transforme plus avant à l'aide des relations de Hecke, qui fournissent alors

(26) $$\Big(\frac{q}{2\pi}\Big)^\delta H(\beta) L(f, \beta+it) L(f, \beta-it) = \sum_{d \geq 1} \frac{\varepsilon(d)}{d} \sum_{n \geq 1} \frac{\lambda_f(n)}{\sqrt{n}} \eta_t(n) V\Big(\frac{2\pi n d^2}{q}\Big).$$

Nous allons exploiter cette formule pour évaluer la somme

(27) $$S := \sum_{f \in S_2(q)^+}^h |L(f, \beta+it) M(f, \beta+it)|^2$$

pour un "mollifier" donné assez arbitrairement, pour commencer, par

$$M(f, \beta+it) = \sum_{m \leq M} \lambda_f(m) x_m m^{-1/2}$$

les $x_m$ étant des nombres complexes quelconques à ceci près que l'on suppose

(28)           Les $x_m$ sont supportés sur les entiers sans facteurs carrés.

(29)           Il existe $A > 0$ tel que $x_m \ll \tau(n)^A (\log n)^A$.

Plus tard nous devrons faire un choix spécifique de $x_m$ pour établir le théorème 1.3, mais le cas général sera utilisé, en vue d'une optimisation, pour étudier la valeur spéciale en $1/2$. Précisons dores et déjà que l'objectif est de démontrer que $S = O(1)$, ou plus exactement que $S = O((|t|+1)^B)$ pour une constante $B > 0$, afin de justifier l'oubli de divers termes d'erreur dans la suite. L'aspect $t$ n'est pas important et la valeur de $B$ pourra évoluer d'une ligne à l'autre.

Il doit être clairement perçu que $M(f, \beta+it)$ dépend holomorphiquement de la variable complexe $\beta + it$ c'est à dire que les $x_i$ en dépendent. Quand le contexte le permettra, on écrira quelquefois $M(f) := M(f, \beta + it)$.

Les relations de Hecke permettent d'écrire cette fois

$$|M(f)|^2 = \sum_b \varepsilon(b) b^{-1} \sum_{m_1, m_2} \lambda_f(m_1 m_2)(m_1 m_2)^{-1/2} x_{bm_1} \overline{x_{bm_2}}$$

(où $x_n = 0$ pour $n > M$ par convention.)

L'emploi de (26) pour exprimer $|L(f, \beta+it)|^2$ fournit donc

$$\begin{aligned}
\left(\frac{q}{2\pi}\right)^\delta H(\beta) S &= \sum_b \frac{\varepsilon(b)}{b} \sum_{m_1, m_2} \sum_n \frac{x_{bm_1} \overline{x_{bm_2}}}{\sqrt{m_1 m_2 n}} \eta_t(n) \sum_d \frac{\varepsilon(d)}{d} V\left(\frac{2\pi d^2 n}{q}\right) \\
&\quad \times \sum_{f \in S_2(q)^+}^h \lambda_f(m_1 m_2) \lambda_f(n). \\
&= \sum_b \frac{\varepsilon(b)}{b} \sum_{m_1, m_2} \sum_n \frac{x_{bm_1} \overline{x_{bm_2}}}{\sqrt{m_1 m_2 n}} \eta_t(n) U\left(\frac{2\pi m_1 m_2}{q}\right) \\
&\quad \times \sum_{f \in S_2(q)^+}^h \lambda_f(m_1 m_2) \lambda_f(n).
\end{aligned}$$

La formule de Petersson (5) s'applique à la somme intérieure sur $f$:

$$\sum_{f \in S_2(q)^+}^h \lambda_f(m_1 m_2) \lambda_f(n) = \delta_{m_1 m_2, n} - \frac{1}{2} \sum_{c \geq 1} (cq)^{-1} S(m_1 m_2, n; cq) J_1\left(\frac{4\pi \sqrt{m_1 m_2 n}}{cq}\right)$$

et on peut utiliser (6) pour estimer la série

$$\sum_{f \in S_2(q)^+}^h \lambda_f(m_1 m_2) \lambda_f(n) = \delta_{m_1 m_2, n} + O((m_1 m_2 n)^{1/2+\varepsilon}(m_1 m_2, n, q)^{1/2} q^{-3/2}).$$

Le terme d'erreur est donc

$$\ll \frac{1}{q^{3/2}} \sum_{b \leq M} b^{-1} |\sum_{bm \leq M} x_{bm}|^2 \times q^{1+\varepsilon}(|t|+1)^B e^{-\pi |t|}$$

pour un certain $\varepsilon > 0$ en utilisant (14) dans le lemme 3.1.

La condition (29) permet de conclure que le terme d'erreur peut être négligé si $M = q^{1/4-\epsilon}$, avec $\epsilon > 0$.

On étudie donc désormais $S_1$ définie par

$$(30) \qquad \left(\frac{q}{2\pi}\right)^\delta H(\beta) S_1 = \sum_b \frac{\varepsilon(b)}{b} \sum_{m_1, m_2} \frac{x_{bm_1}\overline{x_{bm_2}}}{m_1 m_2} \eta_t(m_1 m_2) U\left(\frac{2\pi m_1 m_2}{q}\right)$$

que l'on considère comme une forme quadratique en les $x_m$ que l'on tente de diagonaliser afin d'en simplifier l'expression, où de l'optimiser sous une contrainte linéaire.

Pour cela on applique (20) afin de séparer les variables $m_1$ et $m_2$

$$\left(\frac{q}{2\pi}\right)^\delta H(\beta) S_1 = \sum_b \frac{\varepsilon(b)}{b} \sum_a \mu(a) \sum_{m_1, m_2} \frac{x_{abm_1}\overline{x_{abm_2}}}{a^2 m_1 m_2} \eta_t(m_1)\eta_t(m_2)$$
$$\times U\left(\frac{2\pi}{q} m_1 m_2 a^2\right)$$

### 3.4 Estimation quand $\delta \neq 0$

On suppose ici $\delta = 1/\log q$: c'est le cas dans la preuve du théorème de densité.

On insère (16) dans la dernière formule pour $S_1$; le terme d'erreur est alors clairement

$$\ll \frac{M^{1+\epsilon}}{\sqrt{q}}(|t|+1)^B e^{-\pi|t|} \ll (|t|+1)^B e^{-\pi|t|}$$

compte tenu de la contrainte sur $M$ mise à jour lors de l'étape précédente de l'estimation.

Nous pouvons alors dans le terme principal effectuer la diagonalisation désirée: on regroupe $a$ et $b$ en posant $k = ab$, ce qui amène à introduire la fonction arithmétique

$$\nu_\delta(k) = \frac{1}{k} \sum_{ab=k} \varepsilon(b)\mu(a) a^{-(1+2\delta)}$$

et $S_2$ telle que

$$\left(\frac{q}{2\pi}\right)^\delta H(\beta) S_2 = \left(\frac{q}{2\pi}\right)^\delta H(\beta) \zeta_q(1+2\delta) \sum_k \nu_\delta(k) \Big|\sum_m \eta_t(m) m^{-1-\delta} x_{km}\Big|^2$$
$$+ \left(\frac{q}{2\pi}\right)^{-\delta} H(1-\beta) \zeta_q(1-2\delta) \sum_k \nu_{-\delta}(k) \Big|\sum_m \eta_t(m) m^{-1+\delta} x_{km}\Big|^2$$

de sorte que $S = S_2 + O((|t|+1)^B q^{-\varepsilon})$.

Mais, suivant Selberg, on remarque que pour $0 \leq \delta < 1/2$, on a les inégalités

$$\zeta_q(1-2\delta) \leq 0$$

$$H(1-\beta) = |\Gamma(\tfrac{3}{2} - \beta + it) G(1-\beta+it)|^2 > 0$$

ainsi que

$$\nu_{-\delta}(k) = \frac{1}{k} \prod_{\substack{p|k \\ p\neq q}} (1 - p^{-1+2\delta}) \geq 0$$

ce qui permet de négliger par positivité le second terme et de diviser par $H(\beta)$, avec pour conséquence
$$S_2 \leq \zeta_q(1+2\delta)\sum_k \nu_\delta(k)|y_k|^2$$
à l'aide du changement de variable

(31) $$y_k = \sum_{m \leq M} \eta_t(m) m^{-1-\delta} x_{km}.$$

On va maintenant choisir les $x_m$ pour l'application au théorème 1.3. Pour la méthode employée par Selberg, il faut que le "mollifieur" $M(f)$ soit un polynôme de Dirichlet approchant l'inverse de $L(f)$. On prend donc (comparer au lemme 2.1)

(32) $$x_m = \mu(m) m^{-\delta-it} \sum_n \frac{\varepsilon(n)\mu(mn)^2}{n^{1+2\delta+2it}} g_M(mn)$$

où on introduit pour tout $M \geq 1$ la fonction

(33) $$g_M(x) = \begin{cases} 1, \text{ si } x \leq \sqrt{M} \\ \frac{\log(x/M)}{\log(1/\sqrt{M})}, \text{ si } \sqrt{M} \leq x \leq M \\ 0, \text{ si } x > M \end{cases}$$

Les conditions (28) et (29) sont alors immédiates.

**Proposition 3.1** *Pour $x_m$ donné par (32), on a*
$$S \ll (|t|+1)^B$$
*où $B > 0$ est une constante absolue.*

*Preuve.* — Par définition, on a
$$x_{km} = \frac{\mu(k)}{k^{\delta+it}} \times \frac{\mu(m)}{m^{\delta+it}} \sum_n \varepsilon(n)\mu(kmn)^2 g_M(kmn) n^{-(1+2\delta+2it)}$$
donc
$$y_k = \frac{\mu(k)}{k^{\delta+it}} \sum_{m,n} \frac{\varepsilon(n)\mu(kmn)^2 \mu(m)\eta_t(m) n^{-it}}{(mn)^{1+2\delta+it}} g_M(kmn).$$

On considère séparément les intervalles $1 \leq k \leq M^{1/2}$ et $M^{1/2} \leq k \leq M$. Dans le premier cas, par intégration complexe encore une fois, nous obtenons pour $k \leq \sqrt{M}$
$$\mu(k) k^{\delta+it} y_k = \frac{1}{2i\pi} \int_{(2)} L_k(s+1+2\delta+it) \frac{(M/k)^s(M^s-1)}{\log M} s^{-2} ds$$

la série de Dirichlet ad-hoc étant maintenant
$$L_k(s) = \sum_{l \geq 1} \mu(kl)^2 \Big(\sum_{mn=l} \varepsilon(n) n^{-it} \mu(m)\eta_t(m)\Big) l^{-s}.$$

Or on a
$$\sum_{n\geq 1} \varepsilon(n)n^{-it-s} = \zeta_q(s+it) = \prod_{p\neq q}(1-p^{-s-it})^{-1}$$
et
$$\sum_{m\geq 1} \mu(m)\eta_t(m)m^{-s} = \prod_p (1-p^{-s}(p^{it}+p^{-it}))$$

donc en effectuant le produit puis en extrayant les termes premiers à $k$ et sans facteurs carrés (effet de l'insertion de $\mu(kl)^2$ dans une série de Dirichlet), on obtient la formule très simple
$$L_k(s) = (1-\eta_t(q)q^{-s})\zeta_{qk}(s-it)^{-1}$$

(observons en passant que $k \leq M < q$ donc $(k,q) = 1$.)

Notre traitement est copié sur celui de Luo [Lu]. Nous déformons le contour alors en la courbe $C$ constituée de la droite verticale $\text{Re}(s) = 0$ de $-i\infty$ en $-i\kappa$, le demi-cercle de rayon $\kappa$ situé dans le demi-plan gauche, et la droite de nouveau entre $i\kappa$ et $+i\infty$, $\kappa$ étant assez petit pour que $\zeta$ n'ait pas de zéro à droite de $C$ et que de plus, par la théorie classique on ait
$$L_k(s) \ll \log(|\text{Im}(s) - t| + 2)$$

sur $C$.

On pique alors le pôle simple en $s = 0$, ce qui fournit l'égalité
$$\begin{aligned}\mu(k)k^{\delta+it}y_k &= (1-\eta_t(q)q^{-1-2\delta-it})\zeta_{qk}(1+2\delta) \\ &\quad + \frac{1}{2i\pi}\int_C L_k(s+1+2\delta+it)\frac{(M/k)^s(M^s-1)}{\log M}s^{-2}ds\end{aligned}$$

ce qui donne
$$\mu(k)k^{\delta+it}y_k \ll \frac{1}{\log q}$$

puis, pour la partie $k \leq M^{1/2}$ de la somme $S_2$
$$\begin{aligned}\zeta_q(1+2\delta)\sum_{k\leq M^{1/2}}\nu_\delta(k)\mu(k)^2 k^{-2\delta}|y_k|^2 &\ll (\log q)^{-1}\sum_{k\leq M^{1/2}}\mu(k)^2\nu_\delta(k)k^{-2\delta} \\ &\ll 1.\end{aligned}$$

Pour la partie $M^{1/2} < k \leq M$, la méthode est similaire et seule la représentation intégrale change.

$\square$

D'autre part pour $\beta \geq 2$ on a
$$|M(f,\beta+it)L(f,\beta+it)|^2 - 1 = \sum_{n>\sqrt{M}} c_n n^{-\beta-it}$$

pour certains $c_n \ll \tau(n)^B(\log n)^C$, de sorte que pour un tel $\beta$
$$\sum_f^h |M(f)L(f,\beta+it)|^2 - 1 = O(M^{(\beta-1)/2})$$

En prenant $M = q^{1/4-\varepsilon}$ on en déduit alors du lemme précédent, par un argument de convexité légalité

$$\sum_f^h |M(f, \beta + it)L(f, \beta + it)|^2 - 1 = O((1 + |t|)^B q^{-(1/8 - 2\varepsilon)(\beta - 1/2)}),$$

uniformément pour $\beta \geq 1/2 + 1/\log q$.

Cette estimation étant acquise, on voit que compte tenu de la définition de $M(f)$ comme quasi-inverse de $L(f)$, le raisonnement de Selberg ([Sel], lemme 14 et théorème 4) s'applique sans changement et fournit la première partie du théorème 1.3.

## 4 Théorème de non-annulation en moyenne harmonique

Dans cette section, nous démontrons la minoration (3). Il s'agit d'une légère modification de la section précédente, nous aurons cependant besoin d'un autre estimation beaucoup plus facile: on considère à nouveau un "mollifieur" de longueur $M = q^{1/4-\epsilon}$,

$$M(f, 1/2) = \sum_{m \leq M} \frac{x_m}{m^{1/2}} \lambda_f(m),$$

où les $x_m$ $m \leq M$ sont supposés réels, vérifient les conditions (28), (29) ainsi que la condition de normalisation que l'on rencontrera au cours de la preuve de la proposition 4.1.

$$(34) \qquad \sum_{m \leq M} \frac{x_m}{m} = 1.$$

On va montrer la proposition suivante

**Proposition 4.1** *Pour $M \leq q^{1/4-\epsilon}$, et pour tout vecteur $X = (x_m)_{m \leq M}$ vérifiant les conditions (28), (29), (34) on a l'égalité*

$$(35) \qquad \sum_f^h L(f, 1/2) M(f, 1/2) = 1 + O_\varepsilon(q^{-\varepsilon/2}).$$

*D'autre part, soit $X = (x_m)_{m \leq M}$ le vecteur défini dans la Proposition (4.2) ci-dessous, alors on a l'égalité*

$$(36) \qquad \sum_f^h |L(f, 1/2) M(f, 1/2)|^2 = 2(1 + \frac{\log q}{2 \log M})(1 + O(\frac{\log \log M}{\log M})).$$

La minoration (3) résulte alors de cette cette proposition par application de l'inégalité de Cauchy-Schwarz:

$$\begin{aligned}
(1 + O_\varepsilon(q^{-\varepsilon}))^2 \quad &= (\sum_f^h L(f, 1/2) M(f, 1/2))^2 \\
&\leq \sum_{\substack{f \\ L(f,1/2) \neq 0}}^h 1 \times \sum_f^h |L(f, 1/2) M(f, 1/2)|^2 \\
&\leq 2(1 + \frac{\log q}{2 \log M})(1 + O(\frac{\log \log M}{\log M})) \sum_{\substack{f \\ L(f,1/2) \neq 0}}^h 1
\end{aligned}$$

Comme $M = q^{1/4-\epsilon}$, le résultat en découle.

### 4.1 Evaluation du moment d'ordre $1$

Nous montrons l'égalité (35).

Soit
$$J_{1/2}(f) := \frac{1}{2i\pi} \int_{(2)} \Lambda(f,s) G(s) \frac{ds}{s-1/2},$$

on a alors (rappelons que $G(1/2)\Gamma(1) = 1$)

$$(\frac{q}{2\pi})^{1/4} L(f, 1/2) = (1+\epsilon_f) J_{1/2}(f) = (\frac{q}{2\pi})^{1/4}(1+\epsilon_f) \sum_l \frac{\lambda_f(l)}{l^{1/2}} \widetilde{W}(\frac{(2\pi)^{1/2} l}{q^{1/2}})$$

avec
$$\widetilde{W}(y) = \frac{1}{2i\pi} \int_{(2)} \Gamma(s+1) G(s+1/2) y^{-s} \frac{ds}{s}$$

Par déformation d'un contour on a

(37) $$\widetilde{W}(y) = O(y^{-1}), \ y \to +\infty,$$

(38) $$\widetilde{W}(y) = 1 + O(y), \ y \to 0.$$

On a alors

$$\sum_f^h L(f, 1/2) M(f, 1/2) = \sum_m \frac{x_m}{m^{1/2}} \sum_l l^{-1/2} \widetilde{W}(\frac{(2\pi)^{1/2} l}{q^{1/2}}) \times$$

$$\left( \sum_f^h \lambda_f(l) \lambda_f(m) + \sum_f^h \epsilon_f \lambda_f(l) \lambda_f(m) \right)$$

la première somme de la parenthèse précédente vaut (formule de Petersson)

(39) $$\delta_{l,m} + O(\frac{(lm)^{1/2+\epsilon}}{q^{3/2}}).$$

Rappelons que pour $f \in S_2(q)^+, \epsilon_f = q^{1/2} \lambda_f(q)$, et que d'autre part pour tout $l$, $\lambda_f(l) \lambda_f(q) = \lambda_f(lq)$, on a donc l'égalité

$$\sum_f^h \epsilon_f \lambda_f(m) \lambda_f(l) = \sum_f^h \lambda_f(ql) \lambda_f(m).$$

Remarquons que $m \leq M < q$ de sorte que $\delta_{m,lq} = 0$ et on obtient donc par (5)

$$\sum_f^h \epsilon_f \lambda_f(m) \lambda_f(l) = O(q^{1/2} \sum_c \frac{S(m, ql; cq)}{cq} J_1(\frac{4\pi (mlq)^{1/2}}{q})),$$

remarquant alors que $(m,q) = 1$ et que $S(m,ql;q) = -1$ (c'est une somme de Ramanujan) on obtient que

$$\sum_{f}^{h} \epsilon_f \lambda_f(m) \lambda_f(l) = O(\frac{(lm)^{1/2+\epsilon}}{q}) \tag{40}$$

et finalement en utilisant (37), (38), (39) et (40)

$$\sum_{f}^{h} L(f,1/2)M(f,1/2) = \sum_{m} \frac{x_m}{m} + O(\frac{M^{1+\epsilon}}{q^{1/2}}),$$

et d'après la condition de normalisation (34) cela conclut la preuve de (35).

### 4.2 Evaluation du moment d'ordre 2

Nous montrons l'égalité (36). On reprend les notations de la section 3: dans ce cas précis on a $\delta = t = 0$ et alors $H(1/2) = 1$. Pour $M \leq q^{1/4-\epsilon}$ on a l'égalité

$$\sum_{f}^{h} |M(f,1/2)L(f,1/2)|^2 = S_1 + O(q^{-\epsilon/2})$$

où $S_1$ est définie en (30).

Appliquant alors le Lemme 3.1, (17), on obtient l'égalité

$$S_1 = \sum_{b \leq M} \frac{1}{b} \sum_{a} \frac{\mu(a)}{a^2} \sum_{m_1} \sum_{m_2} x_{abm_1} x_{abm_2} \frac{\tau(m_1)\tau(m_2)}{m_1 m_2}$$

$$\times (\frac{\phi(q)}{q} \log(\frac{q}{2\pi a^2 m_1 m_2}) + C_q) + O(\frac{M^{1/2+\epsilon}}{q^{1/4}})$$

Posont $\log(Q^2) = \log(q/2\pi) + \frac{q}{\phi(q)} C_q$, on a

$$S_1 = \frac{\phi(q)}{q} S(X) + O(q^{-\epsilon})$$

où $S(X)$ désigne la forme quadratique

$$S(X) = \sum_{k} \nu(k) \sum_{m_1,m_2} x_{km_1} x_{km_2} \frac{\tau(m_1)\tau(m_2)}{m_1 m_2} (\log(Q^2/m_1 m_2) + 2\psi(k)).$$

avec $\nu(k) := \nu_0(k) = \phi(k)/k^2$ et

$$\psi(k) = \frac{1}{\nu(k)} \sum_{ab=k} \frac{1}{b} \frac{\mu(a)}{a^2} \log(1/a) = \sum_{p|k} \frac{\log p}{p-1}.$$

Il reste maintenant à choisir un vecteur $X = (x_m)_{m \leq M}$ vérifiant les conditions (28), (29), (34) qui minimise la forme quadratique $S(X)$.

Faisant le changement de variable $x'_m = m x_m$, on a l'égalité $S(X) = Q(X')$ où $Q(X')$ est la forme quadratique étudiée dans [IS], nous écrivons alors le Corollary 6 de [IS] sous la forme de la Proposition suivante:

**Proposition 4.2** *Il existe $X = (x_m)_{m \leq M}$ vérifiant les conditions (28),(29), (34) tel que*
$$S(X) = 2(1 + \frac{\log q}{2 \log M})\{1 + O(\frac{\log \log M}{\log M})\}.$$

Faisant ce choix de $X$, on en déduit l'égalité (36).

**Remarque**. — L'optimisation d'une forme quadratique plus complexe est effectuée en détail dans la section 6.

## 5 Théorèmes de densité en moyenne naturelle

Nous considérons maintenant la seconde partie du théorème 1.3. Un raisonnement exactement identique à celui des paragraphes précédents montre qu'il découle d'une majoration convenable pour la somme

$$S^n = \sum_{f \in S_2(q)^+} |L(f, \beta + it)M(f, \beta + it)|^2, \quad \beta = \frac{1}{2} + \frac{1}{\log q}$$

où $M(f)$ est donné par (32). On posera ici, pour simplifier

$$\alpha_f = |L(f, \beta + it)M(f, \beta + it)|^2.$$

Pour insérer le facteur $(f, f)$ nous allons exploiter le lien avec le carré symétrique des formes automorphes.

D'après Shimura [Sh] et Gelbart-Jacquet [GJ], on introduit

$$L(\mathrm{Sym}^2 f, s) = \zeta_q(2s) \sum_{n \geq 1} \lambda_f(n^2).$$

**Théorème 5.1** $L(\mathrm{Sym}^2 f)$ *est une fonction entière et*

(41) $$\omega_f := L(\mathrm{Sym}^2 f, 1) = \frac{8\pi^3 (f, f)}{q}.$$

*De plus, la fonction L complétée*

$$\Lambda(\mathrm{Sym}^2 f, s) = \pi^{-3s/2} (q^2)^{s/2} \Gamma(\frac{s+1}{2})^2 \Gamma(\frac{s}{2} + 1) L(\mathrm{Sym}^2 f, s)$$

*possède l'équation fonctionnelle*

$$\Lambda(\mathrm{Sym}^2 f, s) = \Lambda(\mathrm{Sym}^2 f, 1 - s).$$

**Remarque**. — On utilise ici - incidemment - que $q$ est premier pour s'assurer que $f$ ne peut être obtenue par changement de base d'un caractère de Hecke d'un corps quadratique, auquel cas $L(\mathrm{Sym}^2 f)$ aurait un pôle en $s = 1$.

La formule (41) va nous permettre d'insérer le facteur $(f, f)$.

Pour cela on va d'abord remplacer la valeur exacte $L(\mathrm{Sym}^2 f, 1)$ par une somme partielle de la série de Dirichlet. Parce que l'on est au bord de la bande critique, on pourra se contenter (en moyenne) d'une somme très courte.

Il faut pour cela la proposition suivante (une forme faible d'une inégalité de grand crible pour les carrés symétriques.)

**Proposition 5.1** *Soit $N > q^9$. Alors pour toute famille de nombres complexes $(a_n)_{n \leq N}$, on a*
$$\sum_{f \in S_2(q)^+} \Big| \sum_{n \leq N} a_n \lambda_f(n) \Big|^2 \ll N(\log N)^{16} \sum_n |a_n|^2.$$

C'est une variante des résultats de [D-K], plus simple en ce que le niveau est fixe et premier, donc le carré symétrique détermine $f$ sans ambiguïté.

Notons $\rho_f(n)$ les coefficients de la série de Dirichlet $L(\text{Sym}^2 f)$ et posons pour tout $x$, $y$, $1 < x < y$
$$\omega_f(y) = \sum_{n \leq y} \rho_f(n) n^{-1}$$

et
$$\omega_f(x,y) = \sum_{x < n \leq y} \rho_f(n) n^{-1}.$$

L'équation fonctionnelle de $L(\text{Sym}^2 f)$ permet d'obtenir alors par transformée de Mellin
$$\omega_f = \omega_f(y) + O(q^{3/2} y^{-1}).$$

On prend $y = q^{5/2}$, $x = q^\epsilon$ et alors on écrit
$$\omega_f = \omega_f(x) + \omega_f(x,y) + O(q^{-1})$$

et comme, selon la proposition 3.1

(42) $$\sum_f^h \alpha_f \ll (|t|+1)^B$$

il vient
$$\frac{q}{2\pi^2} S^n = \sum_{f \in S_2(q)^+}^h \omega_f(x) \alpha_f + \sum_{f \in S_2(q)^+}^h \omega_f(x,y) \alpha_f + O(q^{-1}(|t|+1)^B)$$
$$= S^n(x) + S^n(x,y) + O(q^{-1}(|t|+1)^B)$$

disons. On va traiter $S^n(x)$ et $S^n(x,y)$ séparément. Les deux énoncés suivants terminent alors la preuve, par la même méthode de Selberg que dans la section précédente.

**Lemme 5.2** *On a pour ces valeurs de $x$ et $y$, pour un $\varepsilon' > 0$ dépendant de $\varepsilon$ assez petit*
$$S^n(x,y) \ll (|t|+1)^B q^{-\varepsilon'}.$$

**Lemme 5.3** *On a pour cette valeur de $x$, si $\varepsilon$ est assez petit,*
$$S^n(x) \ll (|t|+1)^B.$$

*Preuve.* — (du lemme 5.2).

La somme $S^n(x,y)$ se traite par la proposition 5.1 et l'inégalité de Hölder: pour tout entier $r \geq 1$, on a

$$\begin{aligned}
\sum_{f \in S_2(q)^+}^h \omega_f(x,y)\alpha_f &= \sum_{f \in S_2(q)^+} \frac{1}{4\pi(f,f)}\omega_f(x,y)\alpha_f \\
&\leq \Big(\sum_f \Big(\frac{\alpha_f}{4\pi(f,f)}\Big)^{\frac{2r}{2r-1}}\Big)^{1-\frac{1}{2r}} \Big(\sum_f \omega_f(x,y)^{2r}\Big)^{\frac{1}{2r}} \\
&\ll \Big(\sum_f^h \alpha_f\Big)^{1-\frac{1}{2r}} \Big(A \sum_f^h \omega_f(x,y)^{2r}\Big)^{\frac{1}{2r}} \\
&\ll (|t|+1)^B \Big(A \sum_f \omega_f(x,y)^{2r}\Big)^{\frac{1}{2r}}
\end{aligned}$$

où on a encore exploité (42) et

$$A = \max_f \frac{\alpha_f}{4\pi(f,f)} \ll (|t|+1)^B \frac{(q^{1/4}M^{1/2})^{2+\varepsilon}}{q} \ll (|t|+1)^B q^{-\gamma}$$

pour un certain $\gamma > 0$, en estimant $\alpha_f$ par convexité par exemple et le fait que pour $q$ premier $(f,f) \gg q/\log q$ [HL, GHL].

La multiplicativité de $\rho_f$ permet d'affirmer l'existence de réels $c(m,n)$ (dépendant de $r$, $x$, et $y$ mais pas de $f$) et de $D > 0$ tels que

(43) $$\omega_f(x,y)^r = \sum_{x^r < mn \leq y^r} \rho_f(m)c(m,n)$$

et

$$|c(m,n)| \leq (mn)^{-1}\tau(mn)^D$$

la variable $n$ étant de plus supportée sur les entiers vérifiant

$$p \mid n \Rightarrow p^2 \mid n.$$

La borne de Deligne $\lambda_f(n) \leq \tau(n)$ permet alors d'estimer la partie avec $n > z$ de (43) comme étant

$$\ll z^{-1}(\log q)^B$$

donc

$$\omega_f(x,y)^r < \sum_{n \leq z} \Big|\sum_{x^r < mn \leq y^r} \rho_f(m)c(m,n)\Big| + O(z^{-1}(\log q)^B).$$

Un découpage dyadique de la somme en $n$ suivi d'une application de l'inégalité de Cauchy permet alors d'établir qu'il existe un entier $N$ tel que

$$x^r z^{-1} < N < y^r$$

et

$$\sum_{f \in S_2(q)^+} \omega_f(x,y)^{2r} \ll \sum_{f \in S_2(q)^+} \Big|\sum_{N < m \leq 2N} c(m)\rho_f(m)\Big|^2 (\log q)^{B_1} + z^{-2}q(\log q)^{B_2}$$

où les réels $c(m)$ ne dépendent pas de $f$ et vérifient

$$|c(m)| \leq m^{-1}\tau(m)^B.$$

L'inégalité de la proposition 5.1 dévoile alors que

$$\sum_f \omega_f(x,y)^{2r} \ll (|t|+1)^B(1+qz^{-2})(\log q)^B$$

si toutefois $x^r \geq zq^9$. On peut donc prendre $z = q^{1/2}$ et pour $r$ assez grand ($r = [10\epsilon^{-1}]+1$ suffit) cela donne

$$\sum_{f \in S_2(q)^+} \omega_f(x,y)^{2r} \ll (|t|+1)^B(\log q)^B$$

d'où le résultat.

$\square$

## 5.1 Preuve du second lemme

On réécrit $\omega_f(x)$ sous la forme

$$\omega_f(x) = \sum_{dl^2 \leq x} \frac{\varepsilon(l)}{dl^2}\lambda_f(d^2)$$

et

$$S^n(x) := \sum_{f \in S_2(q)^+}^h \omega_f(x)|L(f, \beta+it)M(f, \beta+it)|^2.$$

On a alors

$$H(\beta)(\frac{q}{2\pi})^\delta S^n(x) = \sum_b \frac{\varepsilon(b)}{b} \sum_{m_1,m_2} \frac{x_{bm_1}\overline{x_{bm_2}}}{(m_1m_2)^{1/2}} \sum_n \frac{\eta_t(n)}{n^{1/2}}U(\frac{2\pi n}{q})$$

$$\times \sum_{dl^2 \leq x} \frac{\varepsilon(l)}{dl^2} \sum_f^h \lambda_f(d^2)\lambda_f(m_1m_2)\lambda_f(n)$$

Utilisant alors le lemme 3.1, on en déduit comme dans la section précédente que pour tout $\epsilon_1 > 0$

$$H(\beta)(\frac{q}{2\pi})^\delta S^n(x) = \sum_b \frac{\varepsilon(b)}{b} \sum_{m_1,m_2} \frac{x_{bm_1}\overline{x_{bm_2}}}{m_1m_2} \sum_{r|m_1m_2} \varepsilon(r)r$$

$$\times \sum_{\substack{dl^2 \leq x \\ r|d^2}} \frac{\varepsilon(l)}{(dl)^2}\eta_t(\frac{m_1m_2d^2}{r^2})U(\frac{2\pi m_1m_2d^2}{qr^2})$$

$$+O_\varepsilon((|t|+1)^B e^{-\pi|t|}\frac{(xqM^2)^{1+\epsilon_1}}{q^{3/2}})$$

Le terme d'erreur est donc acceptable pour

$$M \leq q^{1/4-\epsilon_2},$$

pour tout $\epsilon_2 > 0$, comme dans la première partie, si l'on prend effectivement $x = q^\epsilon$ avec $\epsilon$ assez petit (en particulier $\varepsilon(b) = \varepsilon(r) = 1$).

On suppose maintenant que $\delta \neq 0$. Par le lemme 3.1 encore, on voit que pour ces valeurs

$$\begin{aligned}(44) \quad H(\beta)(\frac{q}{2\pi})^\delta S^n(x) &= H(\beta)(\frac{q}{2\pi})^\delta \zeta_q(1+2\delta) S_\delta \\ &\quad + H(1-\beta)(\frac{q}{2\pi})^{-\delta} \zeta_q(1-2\delta) S_{-\delta} \\ &\quad + O((|t|+1)^B e^{-\pi t})\end{aligned}$$

où on a posé

$$(45) \quad S_\delta(x) = \sum_b \frac{1}{b} \sum_{m_1,m_2} \frac{x_{bm_1}\overline{x}_{bm_2}}{(m_1 m_2)^{1+\delta}} \sum_{r|m_1 m_2} f_\delta(r) \eta_t(\frac{m_1 m_2}{r}) g_{x,\delta}(r)$$

avec

$$(46) \quad f_\delta(r) = \sum_{r_1 r_2 = r} \mu(r_1) r_2^{1+2\delta}$$

et

$$(47) \quad g_{x,\delta}(r) = \sum_{\substack{dl^2 \leq x \\ r|d^2}} \frac{1}{d^{2+2\delta} l^2} \eta_t(\frac{d^2}{r})$$

Pour montrer (44) et (45) où a utilisé la relation pour $u$ et $v$ fixés

$$\sum_{r|(u,v)} r^{1+2\delta} \eta_t(\frac{uv}{r^2}) = \sum_{r|(u,v)} f_\delta(r) \eta_t(\frac{u}{r}) \eta_t(\frac{v}{r}).$$

**Calcul de $g_{x,\delta}(r)$**

On définit les deux fonctions multiplicatives

$$N(r) = \prod_{p|r} p, \quad M(r) = \prod_{p\|r} p.$$

Pour $r|m_1 m_2$ avec $m_1$ et $m_2$ sans facteurs carrés, on remarque que

$$r|d^2 \iff N(r)|d,$$

et

$$r = M(r)(\frac{N(r)}{M(r)})^2 = \frac{N(r)^2}{M(r)}.$$

Cela permet d'écrire et donc

$$g_{x,\delta}(r) = \frac{1}{N(r)^{2+2\delta}} \sum_{dl^2 \leq x/N(r)} \frac{1}{d^{2+2\delta}l^2} \eta_t(M(r)d^2)$$

$$= g_\delta(r) + O\left(\frac{\log x}{N(r)^{1+2\delta}x^{1/2}}\right)$$

où $g_\delta$ désigne la série étendue à tout $d, l \geq 1$. Le terme d'erreur de l'équation précédente peut être inséré dans (45) sans dommage dès que $\log x \gg \log q$. On y remplace donc $g_{x,\delta}$ par $g_\delta$.

**Lemme 5.4** *On a*

(48) $$g_\delta(r) = g_\delta(1) N(r)^{-2(1+\delta)} \prod_{p||r} \frac{\eta_t(p)}{1 + p^{-2(1+\delta)}}$$

*et*

(49) $$g_\delta(1) = \frac{\zeta(2)\zeta(2+2\delta)|\zeta(2+2\delta+2it)|^2}{\zeta(4+4\delta)} \gg 1$$

*uniformément pour $t \in \mathbf{R}$, $|\delta| \leq \frac{1}{4}$.*

*Preuve.* — On a comme ci-dessus

$$g_\delta(r) = N(r)^{-2(1+\delta)} \zeta(2) g(r, 1+\delta)$$

où

$$g(r,s) = \sum_{d \geq 1} \eta_t(M(r)d^2) d^{-2s}.$$

On a alors

$$g(r,s) = \left(\sum_{\substack{d, \\ (d,M(r))=1}} \eta_t(d^2) d^{-2s}\right) \left(\sum_{d|M(r)^\infty} \eta_t(M(r)d^2) d^{-2s}\right)$$

$$= \tilde{L}(Sym^2\eta_t, 2s) \prod_{p||r} \tilde{L}_p(Sym^2\eta_t, 2s)^{-1}$$

$$\times \prod_{p||r} \sum_{k \geq 0} \eta_t(p^{2k+1}) p^{-2ks}$$

où on a noté $\tilde{L}(Sym^2\eta_t)$ le "faux" carré symétrique de $\eta_t$ (cf. (24):

$$\tilde{L}(Sym^2\eta_t, s) = \sum_{d \geq 1} \frac{\eta_t(d^2)}{d^s} = \frac{\zeta(s-2it)\zeta(s)\zeta(s+2it)}{\zeta(2s)}$$

et $\tilde{L}_p(Sym^2\eta_t)$ son $p$-facteur eulérien. Or

$$\sum_{k \geq 0} \eta_t(p^{2k+1}) p^{-2ks} = \sum_{k \geq 0} (\eta_t(p)\eta_t(p^{2k}) - \eta_t(p^{2k-1})) p^{-2ks}$$

donc
$$(1+p^{-s})\sum_{k\geq 0}\eta_t(p^{2k+1})p^{-2ks} = \eta_t(p)\tilde{L}_p(Sym^2\eta_t, 2s)$$

ce qui donne
$$g(r,s) = \tilde{L}(Sym^2\eta_t, 2s)\prod_{p||r}\frac{\eta_t(p)}{1+p^{-2s}}$$

et en reportant ci-dessus, le lemme.

□

**Diagonalisation de $S_\delta$**

Pour tout couple $m_1, m_2$ apparaissant dans (45), on décompose $m_1m_2 = m_1'm_2'(m_1,m_2)^2$ et tout diviseur $r$ de $m_1m_2$ se décompose de manière unique sous la forme $r = r_1r_2r_3$ avec $r_1|m_1'$, $r_2|m_2'$ $r_3|(m_1,m_2)^2$ observons encore que $m_1', m_2', (m_1,m_2)$ sont premiers entre eux deux à deux ($m_1, m_2$ sont sans facteurs carrés) si bien qu'en utilisant la multiplicativité des fonctions $\eta_t, f_\delta, \frac{g_\delta}{g_\delta(1)}$, on obtient

$$S_\delta(x) = g_\delta(1)\sum_b \frac{1}{b}\sum_{\substack{m_1,m_2\\(m_1,m_2)=1}}\sum_{m_3} h_\delta(m_1)h_\delta(m_2)h_\delta(m_3^2)x_{bm_3m_1}\overline{x}_{bm_3m_2}$$

$$= \sum_b \frac{1}{b}\sum_a \mu(a)h_\delta(a)^2\sum_{m_3}h(m_3^2)\sum_{m_1,m_2}h_\delta(m_1)h_\delta(m_2)x_{bm_3am_1}\overline{x}_{bm_3am_2}$$

(50) $$= \sum_k \tilde{\nu}(k)\Big|\sum_m h_\delta(m)x_{km}\Big|^2$$

avec
(51) $$h_\delta(m) = \frac{1}{m^{1+\delta}}\sum_{r_1m_1=m}f_\delta(r_1)\frac{g_\delta(r_1)}{g_\delta(1)}\eta_t(m_1)$$

et

(52) $$\tilde{\nu}_\delta(k) = \sum_{abr_3=k}\frac{1}{b}\mu(a)h_\delta(a)^2 h_\delta(r_3^2)$$

Notons les égalités
$$h_\delta(p) = \frac{\eta_t(p)}{p^{1+\delta}}(1+\frac{1-p^{-1-2\delta}}{p(1+p^{-2-2\delta})})$$

$$h_\delta(p^2) = \frac{1}{p^{2+2\delta}}(p^{2\delta}-1-\frac{1}{p}+\eta_t^2(p)(1+\frac{1-p^{-1-2\delta}}{p(1+p^{-2-2\delta})}))$$

et
$$\tilde{\nu}_\delta(p) = \frac{1}{p} - h_\delta^2(p) + h_\delta(p^2)$$

On peut alors vérifier que si $q$ est assez grand, pour tout $t \in \mathbf{R}$ et tout $\delta$ tel que $|\delta| \leq 1/\log q$, on a $\tilde{\nu}_\delta(k) > 0$. Donc, en vertu de cette remarque et de (49), on voit qu'à un

terme admissible près $S_{-\delta}$ est positive ou nulle. On conclut donc comme précédemment en utilisant la négativité de $\zeta_q(1-2\delta)$ que

$$S^n(x) \leq g_\delta(1)\zeta_q(1+2\delta)\sum_k \tilde{\nu}_\delta(k)\left|\sum_r \frac{f_\delta(r)g_\delta(r)/g_\delta(1)}{r^{1+\delta}}\sum_m \frac{x_{krm}\eta_t(m)}{m^{1+\delta}}\right|^2$$

raisonnant comme dans la section précédente on en déduit en utilisant le fait que pour $r$ sans facteurs carrés $f_\delta(r)g_\delta(r)/g_\delta(1) \ll \tau(r)/r^2$ que

$$S^n(x) \ll (|t|+1)^B$$

et donc que
(53) $$S^n \ll q(|t|+1)^B$$

comme désiré.

## 6 Théorème de non-annulation en moyenne naturelle

Dans cette section, nous prouvons la minoration (4) qui est une modification des arguments de la section 4. Dans notre cas on impose au mollifieur la normalisation suivante

(54) $$\sum_{m \leq M} \frac{x_m d_{-1}(m)}{m} = 1, \ avec \ d_{-1}(m) := \sum_{r|m} r^{-1}.$$

On va montrer la proposition

**Proposition 6.1** *Pour $M \leq q^{1/4-\epsilon}$, et pour tout vecteur $X = (x_m)_{m \leq M}$ vérifiant les conditions (28), (29), (54) on a l'égalité*

(55) $$\sum_f^h \omega_f L(f,1/2)M(f,1/2) = \zeta(2)^2 + O(q^{-\epsilon}).$$

*D'autre part, soit $X = (x_m)_{m \leq M}$ le vecteur défini dans par (67) ci-dessous, alors on a l'égalité*

(56) $$\sum_f^h \omega_f |L(f,1/2)M(f,1/2)|^2 = 2\zeta(2)^3(1 + \frac{\log q}{2\log M})(1 + O(\frac{\log\log M}{\log M})).$$

On conclut comme dans la section 4 et considérant l'égalité (41) et en remarquant que $|S_2(q)^+| = \frac{q}{12} + O(1)$.

En raisonnant comme dans la section 5, on voit que pour prouver (55) et (56) il suffit de considérer majorer les quantités

$$\sum_f^h \omega(x)L(f,1/2)M(f,1/2)$$

et

$$\sum_f^h \omega(x)|L(f,1/2)M(f,1/2)|^2,$$

avec $x = q^\varepsilon$, les autres termes fournissant une erreur en $q^{-\varepsilon'}$

### 6.1 Evaluation du moment d'ordre 1

$$\sum_f^h \omega_f(x)L(f,1/2)M(f,1/2) = \sum_m \frac{x_m}{m^{1/2}} \sum_{dl^2 \leq x} \frac{1}{dl^2} \sum_{rr'=d^2} r^{-1/2} \sum_l l^{-1/2}\widetilde{W}(\frac{(2\pi)^{1/2}rl}{q^{1/2}}) \times$$

$$\left(\sum_f^h \lambda_f(r'l)\lambda_f(m) + \sum_f^h \epsilon_f \lambda_f(r'l)\lambda_f(m)\right)$$

Par un raisonnement analogue à celui de la section précédente (cf. (39), (40)) la paranthèse vaut

(57) $$\delta_{r'l,m} + O(\frac{(d^2lm)^{1/2+\epsilon}}{q}).$$

et finalement en utilisant (57), (38) et (37) on obtient

$$\sum_f^h \omega_f(x)L(f,1/2)M(f,1/2) = \sum_{dl^2 \leq x} \frac{1}{d^2l^2} \sum_m \frac{x_m}{m} \sum_{r|(d^2,m)} r + O(\frac{x^2M^{1+\epsilon}}{q^{1/2}}).$$

On remarque que $(d^2,m) = (d,m)$ car $m$ est sans facteur carrés et on obtient alors

$$\sum_f^h \omega_f(x)L(f,1/2)M(f,1/2) = \sum_m \frac{x_m}{m} \sum_{r|m} r^{-1} \sum_{dl^2 \leq x/r} \frac{1}{d^2l^2}$$

$$= \zeta(2)^2 \sum_m \frac{x_m d_{-1}(m)}{m} + O(\frac{(\log q)^B}{x^{1/2}}) + O(\frac{x^2M^{1+\epsilon}}{q^{1/2}}).$$

d'après (54) cela conclut la preuve de (55).

### 6.2 Evaluation du moment d'ordre 2

On se trouve dans le cas $\delta = t = 0$ et les égalités (46), (48), (49) et (52) deviennent

$$f_0(r) = \phi(r)$$

(58) $$g_0(r) := g(r) = \frac{\zeta(2)^4}{\zeta(4)N(r)^2} \prod_{p||r} \frac{2}{1+p^{-2}}$$

(59) $$h_0(p) := h(p) = \frac{2}{p}(1 + \frac{1-p^{-1}}{p(1+p^{-2})})$$

(60) $$\tilde{\nu}(k) := \tilde{\nu}_0(k) = \frac{1}{k} \prod_{p|k} D(1/p), \text{ avec } D(x) = \frac{(1-x^2)^3}{(1+x^2)^2}$$

De manière analogue à (50), on trouve pour $M \leq q^{1/4-\epsilon}$, l'égalité

$$\sum_f^h \omega_f(x)|M(f,1/2)L(f,1/2)|^2 = S_1^n + O(q^{-\epsilon/2})$$

avec

$$S_1^n = \frac{\zeta(2)^4}{\zeta(4)} \sum_{b,a,r_3} \frac{1}{b}\mu(a)h(a^2)h(r_3^2)$$

$$\times \sum_{m_1}\sum_{m_2} x_{abr_3m_1}x_{abr_3m_2}h(m_1)h(m_2)$$

$$\times (\frac{\phi(q)}{q}\log(\frac{q}{2\pi a^2 r_3^2 m_1 m_2}) + C_q) + O(\frac{M^{1/2+\epsilon}}{q^{1/4}}).$$

On peut encore écrire ceci sous la forme

(61) $$S_1^n = \frac{\phi(q)}{q}\frac{\zeta(2)^4}{\zeta(4)}S^n(X) + O(q^{-\epsilon})$$

où $S^n(X)$ désigne la forme quadratique

$$S^n(X) = \sum_k \tilde{\nu}(k) \sum_{m_1,m_2} x_{km_1}x_{km_2}h(m_1)h(m_2)(\log(Q^2/m_1 m_2) + 2\tilde{\psi}(k)).$$

avec

$$\tilde{\psi}(k) = \frac{1}{\tilde{\nu}(k)} \sum_{abr_3=k} \frac{1}{b}\mu(a)h(a)^2 h(r_3^2)\log(1/ar_3)$$

$$= \sum_{p|k} \log p \frac{1/p - \tilde{\nu}(p)}{\tilde{\nu}(p)} = O(\log\log k).$$

Il reste maintenant à choisir un vecteur $X = (x_m)_{m \leq M}$ vérifiant les conditions (28),(29), (54) qui minimise la forme quadratique $S^n(X)$.

### 6.3 Diagonalisation de $S^n(X)$

Dans cette partie les lettres $k, k', m$ désignerons exclusivement des entiers sans facteurs carrés. Nous suivons de près les argument de [IS]: posons

$$y_k = \sum_{m \leq M/k} h(m)x_{km}$$

de sorte que si on note $h^{(-1)}(m)$ l'inverse de convolution de $h(m)$ ($h^{(-1)}(p) = -h(p)$) on a

(62) $$x_m = \sum_{k \leq M/m} h^{(-1)}(k)y_{km}$$

et la condition (54) est équivalente à

(63) $$\sum_k y_k(h^{(-1)} * \frac{d_{-1}}{Id})(k) = 1.$$

On pose dans la suite $j(k) := (h^{(-1)} * \frac{d_{-1}}{Id})(k)$ et on a

$$j(p) = \frac{-1}{p} E(1/p), \ avec \ E(x) = \frac{(1-x)(1+x)^2}{1+x^2}.$$

On a alors

$$\sum_k \tilde{\nu}(k) \sum_{m_1,m_2} x_{km_1} x_{km_2} h(m_1) h(m_2) \log(\frac{Q^2}{m_1 m_2}) = 2 \sum_k \tilde{\nu}(k) y_k y'_k$$

avec

$$\begin{aligned} y'_k &= \sum_m x_{km} h(m)(\log Q - \log m) \\ &= (\log Q) y_k - \sum_{k'} \Lambda(k') h(k') y_{kk'} \end{aligned}$$

de sorte que

$$S^n(X) = R(Y) - R'(Y)$$

avec
(64)
$$R(Y) = 2 \sum_{k \leq M} \tilde{\nu}(k)(\log Q + \tilde{\psi}(k)) y_k^2$$

et
(65)
$$R'(Y) = 2 \sum_{kk' \leq M} \tilde{\nu}(k) \Lambda(k') h(k') y_k y_{kk'}$$

Pour minimiser $R(Y)$ il est clair que le choix optimal est donné par:

(66)
$$y_k = \frac{\mu^2(k) j(k)}{\tilde{\nu}(k)}/J, \ avec \ J = \sum_k \frac{\mu^2(k) j^2(k)}{\tilde{\nu}(k)}$$

de sorte que la condition (63) est vérifiée.

**Evaluation de $J$**

On a, pour $k$ sans facteurs carrés,

$$\frac{j^2(k)}{\tilde{\nu}(k)} = \frac{1}{k} \prod_{p|k} \frac{1+1/p}{1-1/p}$$

On voit donc que

$$\sum_k \mu^2(k) \frac{j^2(k)}{\tilde{\nu}(k) k^s} = \prod_p (1-1/p)^{-1}(1 - \frac{1}{p} + \frac{1}{p^{1+s}} + \frac{1}{p^{2+s}})$$

et ainsi

$$\begin{aligned} J = R_1 \log M + O(1), \ avec \ R_1 &= res_{s=0} \prod_p (1-1/p)^{-1}(1 - \frac{1}{p} + \frac{1}{p^{1+s}} + \frac{1}{p^{2+s}}) \\ &= \prod_p (1 + \frac{1}{p^2}) = \zeta(2)/\zeta(4). \end{aligned}$$

Alors par (62), on a

$$
\begin{aligned}
x_m &= \frac{\mu^2(m)j(m)}{J\tilde{\nu}(m)} \sum_{\substack{k \leq M/m \\ (k,m)=1}} \frac{\mu^2 j \times h^{(-1)}}{\tilde{\nu}}(k) \\
&= \frac{\mu^2(m)j(m)}{J\tilde{\nu}(m)} \sum_{\substack{k \leq M/m \\ (k,m)=1}} \mu^2(k) \prod_{p|k}(\frac{2}{p} + O(\frac{1}{p^2})) \\
&\ll \tau(m)\log^2 m
\end{aligned}
\tag{67}
$$

Ainsi les conditions (28), (29) et (54) sont vérifiées

**Evaluation de $R'(Y)$**

$$
\begin{aligned}
R'(Y) &= \frac{2}{J^2} \sum_{k \leq M} \mu^2(k)\tilde{\nu}(k)\frac{j(k)}{\tilde{\nu}(k)} \sum_{\substack{p \leq M/k \\ (p,k)=1}} \frac{j(kp)h(p)}{\tilde{\nu}(kp)} \log p \\
&= \frac{2}{J^2} \sum_{k \leq M} \mu^2(k)\frac{j^2(k)}{\tilde{\nu}(k)} \sum_{\substack{p \leq M/k \\ (p,k)=1}} \frac{j(p)h(p)}{\tilde{\nu}(p)} \log p
\end{aligned}
$$

Or, on a $\frac{j(p)h(p)}{\tilde{\nu}(p)} = -\frac{2}{p} + O(\frac{1}{p^2})$, et donc

$$
\sum_{\substack{p \leq M/k \\ (p,k)=1}} \frac{j(p)h(p)}{\tilde{\nu}(p)} \log p = -2\log(M/k) + O(\log\log k)
$$

Or on a l'égalité

$$
\sum_{k \leq M} \mu^2(k)\frac{j^2(k)}{\tilde{\nu}(k)} \log k = \frac{\zeta(2)}{2\zeta(4)} \log^2 M + O((\log\log M)\log M)
$$

et ainsi

$$
-R'(Y) = \frac{2\zeta(4)}{\zeta(2)} + O(\frac{\log\log M}{\log M})
$$

**Evaluation de $R(Y)$**

On trouve de même que

$$
R(Y) = 2\frac{\zeta(4)}{\zeta(2)}\frac{\log Q}{\log M}\left(1 + O(\frac{\log\log M}{\log M})\right),
$$

se qui donne finalement l'égalité

$$
S^n(X) = \frac{2\zeta(4)}{\zeta(2)}(1 + \frac{\log q}{2\log M})(1 + O(\frac{\log\log q}{\log q}))
\tag{68}
$$

Cette dernière égalité conjuguée avec (61) et (41), conclut la preuve de (56)

# 7 Le premier théorème de densité

La méthode de démonstration du théorème 1.1 est basée sur la méthode des "polynômes détecteurs de zéros", due à Montgomery (voir [Mon], [Bo]). Nous allons suivre de très près l'exposition de Bombieri, [Bo], pages 76 à 82.

## 7.1 Variantes du grand-crible

Rappelons que la notation $n \sim N$ signifie $N < n \leq 2N$. La Proposition 2.1 entraîne aussitôt par positivité

$$\sum_{f \in S_2(q)^+}^h \Big| \sum_{n \sim N} \beta_n \lambda_f(n) \Big|^2 \ll (1 + N/q) \sum_{n \sim N} |\beta_n|^2$$

(le poids $(f,f)$ provenant de la renormalisation des formes primitives en système orthonormé) et la constante du $\ll$ est absolue.

La première conséquence est une intégration triviale en $t$:

$$(69) \qquad \int_{-T}^{T} \sum_{f \in S_2(q)^+}^h \Big| \sum_{n \sim N} \beta_n \lambda_f(n) n^{-it} \Big|^2 dt \ll T(1 + N/q) \sum_{n \sim N} |\beta_n|^2.$$

Ensuite on remplace l'intégration en $t$ par une somme sur des points bien espacés (ce seront ultérieurement les zéros des fonctions $L$). Ici $\alpha$ est entre $1/2$ et $1$.

**Lemme 7.1** *Supposons donné pour toute $f \in S_2(q)^+$ un ensemble fini $J_f$ de nombres complexes tous inclus dans la région*

$$R(\alpha, T) = \{ z = \sigma + it \in \mathbf{C} \mid \alpha \leq \sigma \leq 1, |t| \leq T \}$$

*tels que de plus pour $\rho \neq \rho'$ dans $J_f$ on ait*

$$|\text{Im}(\rho - \rho')| \gg 1.$$

*Alors on a*

$$\sum_{f \in S_2(q)^+}^h \sum_{\rho \in J_f} \Big| \sum_{n \sim N} \beta_n \lambda_f(n) n^{-\rho} \Big|^2 \ll T \log N \sum_{n \sim N} |\beta_n|^2 n^{-2\alpha}(1 + n/q)$$

*(toujours pour toute famille de nombres complexes $\beta_n$.)*

La démonstration est la même que celle du théorème 7.3 de [Mon] (on utilise un lemme de Gallagher).

**Lemme 7.2** *Soient $(\beta_n)_{n \geq 1}$ une suite de nombres complexe dans $\ell_1$, $T \geq 1$ un réel. Alors on a*

$$\sum_{f \in S_2(q)^+}^h \int_{-T}^{T} \Big| \sum_{n \geq 1} \beta_n \lambda_f(n) n^{it} \Big|^2 dt \ll T \sum_{n \geq 1} (1 + n/q)|\beta_n|^2$$

Pour la preuve, voir [Bo], page 30. Noter que l'on perd un facteur $T$ par rapport à sa version, faute de disposer du grand crible dans des petits intervalles $M \leq n \leq M + N$.

## 7.2 Lemmes utilitaires

Nous nous intéressons ici à des séries de Dirichlet de la forme

$$\sum_{l,m\geq 1} c(l,m)(lm^2)^{-s},$$

catégorie qui inclus, comme on l'a vu, les inverses et carrés des fonctions $L$ de nos formes primitives.

L'objectif est de voir qu'en certaines circonstances, ces séries se comportent comme si elles étaient réduites à la simple expression

$$\sum_{l\geq 1} c(l,1)l^{-s}$$

(en particulier dans certaines estimations.)

**Lemme 7.3** *Soit* $(c(l,m))_{l,m\geq 1}$ *des nombres complexes et* $M$ *un réel tel que* $m \geq M$ *implique* $c(l,m) = 0$. *Alors pour tout* $s \in \mathbf{C}$

$$\Big| \sum_{l,m\geq 1} c(l,m)(lm^2)^{-s} \Big|^2 \ll (\log M) \sum_{m\leq M} m^{1-4\mathrm{Re}(s)} \Big| \sum_{l\geq 1} c(l,m)l^{-s} \Big|^2.$$

*Preuve.* — C'est clair par Cauchy-Schwarz, en écrivant la somme ($\kappa$ est n'importe quel réel)

$$\Big| \sum_{l,m} c(l,m)(lm^2)^{-s} \Big|^2 = \Big| \sum_{m\leq M} m^{\kappa-s} \sum_{l} m^{-\kappa-s} c(l,m) l^{-s} \Big|^2$$

$$\leq \Big( \sum_{m\leq M} m^{2(\kappa-\mathrm{Re}(s))} \Big) \Big( \sum_{m\leq M} \Big| \sum_{l} m^{-\kappa-s} c(l,m) l^{-s} \Big|^2 \Big)$$

et en choisissant $\kappa$ de sorte que $\mathrm{Re}(s) - \kappa = 1/2$.

$\square$

Le second lemme est une transformation simple basée sur la multiplicativité des coefficients de $L(f)$.

**Lemme 7.4** *Soient* $L$, $M$, $N$ *des réels avec* $LM = N$ *et* $f \in S_2(q)^+$. *Pour toute famille* $(c(l,m))_{l\sim L, m\sim M}$ *de nombres complexes, soit*

$$\tilde{c}(n,d) = \sum_{\substack{lm=n \\ ld\sim L, md\sim M}} \varepsilon(d) c(ld, md).$$

*Alors*

$$\sum_{n\sim N} \Big( \sum_{\substack{lm=n \\ l\sim L, m\sim M}} c(l,m) \lambda_f(l) \lambda_f(m) \Big) n^{-s} = \sum_{nd^2\sim N} \tilde{c}(n,d) \lambda_f(n) (nd^2)^{-s}.$$

*Preuve.* — Par la formule de multiplication des coefficients qui, pour les $\lambda_f(n)$, est

$$\lambda_f(l)\lambda_f(m) = \sum_{d|(l,m)} \varepsilon(d)\lambda_f(lm/d^2).$$

$\square$

On peut mixer les résultats du paragraphe précédent au lemme 7.3. Soient par exemple des complexes $(c(l,m))_{l,m\geq 1}$ vérifiant $|c(l,m)| \leq |c(l,1)|$. On a alors

$$(70) \quad \sum_{f\in S_2(q)^+}^h \sum_{\rho\in J_f} |\sum_{lm^2\leq N} c(l,m)\lambda_f(l)(lm^2)^{-\rho}|^2 \ll \log N(1+N/q)S_\alpha(N)$$

$$\times \sum_{l\leq N} |c(l,1)|^2 l^{-2\alpha}$$

où

$$S_\alpha(N) = \sum_{n\leq N^{1/2}} n^{1-4\alpha}.$$

De même en partant du Lemme 7.2, supposant que $(c(l,1))_{l\geq 1}$ est dans $\ell_2$. Si $\beta > 1/2$ on déduit que

$$(71) \quad \sum_{f\in S_2(q)^+}^h \int_{-T}^T \Big|\sum_{l,m\geq 1} c(l,m)\lambda_f(l)(lm^2)^{-\beta-it}\Big|^2 dt \ll_\beta T\sum_{l\geq 1} |c(l,1)|^2 l^{-2\beta}(1+q/l).$$

### 7.3 L'estimation auxiliaire

L'estimation obtenue dans ce paragraphe, dans la lignée du précédent, servira à traiter l'inverse de $L(f)$ "comme si" il était simplement donné en multipliant les coefficients par $\mu(n)$, comme c'est le cas pour les caractères de Dirichlet (voir la remarque à la fin.) Il s'agit d'une partie purement technique.

On suppose donnés des complexes $c(k,l,m)_{k\sim K,\, l\sim L,\, m\sim M}$ et on veut estimer la somme

$$(72) \quad \sum_{f\in S_2(q)^+}^h \sum_{\rho\in J_f} \Big|\sum_{n\sim N} \Big(\sum_{\substack{klm^2=n \\ k\sim K,\, l\sim L,\, m\sim M}} c(k,l,m)\lambda_f(k)\lambda_f(l)\Big) n^{-\rho}\Big|^2$$

où $J_f$ est un ensemble fini bien espacé, comme dans l'énoncé du lemme 7.1.

On procède en deux étapes.

**Première étape.** On considère seulement $m=1$, c'est à dire une somme

$$\sum_{f\in S_2(q)^+}^h \sum_{\rho\in J_f} \Big|\sum_{n\sim N} \Big(\sum_{\substack{kl=n \\ k\sim K,\, l\sim L}} c(k,l)\lambda_f(k)\lambda_f(l)\Big) n^{-\rho}\Big|^2.$$

On applique le lemme 7.4 pour écrire la somme sur $n$ sous la forme

$$\sum_{n,d} \tilde{c}(n,d)(nd^2)^{-\rho}$$

puis le lemme 7.3 pour séparer $n$ et $d$. On somme ensuite sur $f$ et $\rho$, et pour chaque valeur de $d$ on applique la variante du grand crible du lemme 7.1, ce qui donne l'estimation

$$\ll \log^2 N \sum_{d \leq N^{1/2}} d^{1-4\alpha} T\Big(1 + \frac{Nd^{-2}}{q}\Big) \sum_{n \sim Nd^{-2}} |\tilde{c}(n,d)|^2 n^{-2\alpha}.$$

Nous supposons maintenant que les $c(k,l)$ vérifient de plus

(73) $$\sum_{n \sim Nd^{-2}} |\tilde{c}(n,d)|^2 n^{-2\alpha} \ll \Big(\frac{N}{d^2}\Big)^{1-2\alpha} (\log N)^A$$

pour une certaine constante $A > 0$.

On peut alors continuer l'estimation, obtenant pour la somme

$$\ll (\log N)^{2+A} \sum_{d \leq N^{1/2}} d^{1-4\alpha} T\Big(1 + \frac{Nd^{-2}}{q}\Big)(Nd^{-2})^{1-2\alpha}$$

$$\ll (\log N)^{2+A} T \Big( \sum_{d \leq N^{1/2}} d^{1-4\alpha-2+4\alpha} + \frac{N}{q} \sum_{d \leq N^{1/2}} d^{-1-4\alpha-2+4\alpha} \Big) N^{1-2\alpha}$$

$$\ll (\log N)^{3+A} T(1 + N/q) N^{1-2\alpha}.$$

**Deuxième étape.** On peut maintenant traiter le cas général. Dans la somme originale (72), on utilise le lemme 7.3 pour séparer $m$ cette fois. Pour chaque $m$ on se trouve confronté, en remplaçant $N$ par $Nm^{-2}$, a une somme du type traité précédemment. Admettons que l'hypothèse (73) correspondante est valide uniformément en $m$; on peut alors estimer notre somme par

$$\ll (\log M) \sum_{m \sim M} m^{1-4\alpha} (\log N)^{3+A} T\Big(1 + \frac{Nm^{-2}}{q}\Big)(Nm^{-2})^{1-2\alpha}$$

et par la même sommation que ci-dessus, ceci est (encore)

$$\ll (\log M)^2 (\log N)^{3+A} T(1 + N/q) N^{1-2\alpha}.$$

L'essentiel est que cette borne est la même que celle que l'on aurait obtenue (par le lemme 7.1) pour une somme

$$\sum_f^h \sum_\rho \Big| \sum_{n \sim N} c(n) \lambda_f(n) n^{-\rho} \Big|^2$$

pour des complexes $(c(n))$ vérifiant la condition

$$\sum_{n \sim N} |c(n)|^2 n^{-2\alpha} \ll (\log N)^A N^{1-2\alpha}$$

(si $M$ est au plus une puissance de $N$), condition qui s'interprète intuitivement comme le fait que les $(c(n))$ sont "presque" bornés, en norme $\ell_2$.

### 7.4 Preuve du théorème

Ayant amassé nos munitions, nous pouvons commencer la campagne suivant le plan décrit au début de cette section.

Nous introduisons donc un inverse tronqué de $L(f)$. Soit $X > 1$ un paramètre qui sera choisi ultérieurement et posons

$$M_X(f,s) = \sum_{lm^2 \leq X} \varepsilon(m)\mu(l)|\mu(lm)|\lambda_f(l)(lm^2)^{-s}.$$

Nous effectuons maintenant le produit $M_X(f,s)L(f,s)$. Pour $\mathrm{Re}(s) > 2$ il vient

$$L(f,s)M_X(f,s) = 1 + \sum_{n>X} \left( \sum_{klm^2=n} \varepsilon(m)\mu(l)|\mu(lm)|\lambda_f(l)\lambda_f(k) \right) n^{-s},$$

les termes avec $1 < n \leq X$ disparaissent puisque tous les coefficients de l'inverse sont présents jusque là.

Notons immédiatement

(74) $$b_f(n) = \sum_{\substack{klm^2=n \\ lm^2 \leq X}} \varepsilon(m)\mu(l)|\mu(lm)|\lambda_f(l)\lambda_f(k)$$

et nous pouvons reformuler cela plus simplement en l'égalité

$$L(f,s)M_X(f,s) = 1 + \sum_{n>X} b_f(n)n^{-s}.$$

Notons $\varphi(x) = e^{-x}$, $\hat{\varphi} = \Gamma$ sa transformée de Mellin. Comme dans [Bo], on obtient donc pour tout $\sigma < 0$ cette fois

$$L(f,s)M_X(f,s) = 1 + \sum_{n>X} b_f(n)n^{-s}\varphi(n/Y)$$

$$-\frac{1}{2i\pi}\int_{(-\sigma+1/2)} Y^w L(f,s+w) M_X(f,s+w)\hat{\varphi}(w) dw.$$

où $\sigma = \mathrm{Re}(s)$.

Si on prend $s = \rho = \beta + i\gamma$, un des zéros dans la région concernée par le théorème, on trouve comme Bombieri que $\rho$ satisfait l'une des trois inégalités suivantes:

- 
$$1 \ll \Big| \sum_{n>X} b_f(n)n^{-\rho}\hat{\varphi}(n/Y) \Big|.$$

- 
$$|M_X(f, 1/2 + it_\rho)| > X^{\alpha - 1/2}$$

  pour un certain $t_\rho$ tel que $|t_\rho - \gamma| < (\log qT)^2$.

- $$\int_{\gamma-(\log qT)^2}^{\gamma+(\log qT)^2} |L(f,1/2+it)|dt \gg (\frac{Y}{X})^{\alpha-1/2}.$$

On traite chaque classe de zéro caractérisée par ces inégalités séparément, et pour chacune on partitionne l'ensemble des zéros en sous-ensembles bien espacés (de $(\log qT)^2$) et on somme les estimations obtenues.

Pour la première classe, la condition donne

$$\frac{1}{(f,f)^{1/2}} \ll \frac{1}{(f,f)^{1/2}}\Big|\sum_{n>X} b_f(n)n^{-\rho}\hat\varphi(n/Y)\Big|$$

et après partitions dyadiques, on se trouve avec (à la place des sommes

$$\sum_\chi |\sum_{n\sim N} \mu(n)\chi(n)n^{-\rho}|^2$$

de [Bo]), des sommes

$$\sum_{f\in S_2(q)^+}^h \sum_{\rho\in J_f} \Big|\sum_{n\sim N}\Big(\sum_{\substack{klm^2=n\\k\sim K,l\sim L,m\sim M}} c(k,l,m)\lambda_f(k)\lambda_f(l)\Big)n^{-\rho}\Big|^2$$

avec $c(k,l,m) = \varepsilon(m)\mu(l)|\mu(lm)|\varphi(klm^2/Y)$, si $lm^2 \leq X$, et $c(k,l,m)=0$ sinon, auxquelles s'applique le résultat du paragraphe précédent.

En effet, comme $|c(k,l,m)|\leq 1$, pour tout $m$ on a $\tilde c(n,d)\leq \tau(n)$ et l'hypothèse (73) est aisément vérifiée avec $A=15/2$ par Cauchy-Schwarz et la borne

$$\sum_{n\leq N}\tau(n)^4 \ll n(\log n)^{15},$$

ce qui signifie simplement que l'on obtient exactement la même borne que dans le cas des caractères. Avec les mêmes choix de $X$ et $Y$ que [Bo] ($X=(qT)^{1/(2-\alpha)}$, $Y=X^{3/2}$) et compte tenu de divers $T$ perdus, on trouve

(75) $$\sum_{f\in S_2(q)^+}^h N_1(f,\alpha,T) \ll \frac{1}{q}T^B(qT)^{\frac{3}{2-\alpha}(1-\alpha)}(\log q)^A$$

(l'indice indiquant que c'est la première classe de zéros.)

Pour la seconde, l'argument est plus simple, puisqu'on aura des sommes

$$\sum_{f\in S_2(q)^+}^h \sum_\rho |M_X(f,1/2+it_\rho)|^2$$

qui sont de la forme considérée dans (70) avec $c(l,m)=\varepsilon(m)\mu(l)|\mu(lm)|$. Là encore, on a la même estimation que dans le cas des caractères.

Pour la troisième, un tel zéro satisfait

$$X^{\alpha-1/4}\frac{1}{(f,f)^{1/4}} \ll \int_{\gamma-(\log qT)^2}^{\gamma+(\log qT)^2}\frac{|L(f,1/2+it)|}{(f,f)^{1/4}}dt$$

donc
$$X^{2\alpha-1}\sum_{f\in S_2(q)^+}^{h} N_3(f,\alpha,T) \ll (\log qT)^6 \sum_{f\in S_2(q)^+}^{h} \int_{-T}^{T} |L(f,1/2+it)|^4 dt$$

et le théorème ne requiert plus que de prouver l'analogue du lemme sur le moment d'ordre 4:

**Lemme 7.5** *Il existe $A$ et $B$ absolus tel que*
$$\sum_{f\in S_2(q)^+}^{h} \int_{-T}^{T} |L(f,1/2+it)|^4 dt \ll qT^B (\log q)^A.$$

*Preuve.* — On suit encore Bombieri, [Bo] page 80-82. L'argument s'adapte ici parce qu'on ne s'intéresse qu'à la dépendance en $q$. On utilise la série de Dirichlet donnée par le lemme 2.2 pour le carré de $L(f)$, et l'équation fonctionnelle sous la forme
$$L(f,s) = \psi(f,s) L(1-s,f)$$
avec
$$\psi(s,f) = g(f) q^{1/2-s} (2\pi)^{s-1/2} \frac{\Gamma(3/2-s)}{\Gamma(1/2+s)}.$$

On écrit (avec $Z = qT$)
$$L(f,s)^2 = \sum_{n\geq 1} c_f(n) n^{-s} \varphi(n/Z) - \frac{1}{2i\pi} \int_{(c)} L(f,s+w)^2 Z^w \hat{\varphi}(w) dw$$

on applique l'équation fonctionnelle et l'intégrale est encore
$$\frac{1}{2i\pi}\int_{(c)} \psi(f,s+w)^2 \Big(\sum_{n>Z}()\Big) Z^w \hat{\varphi}(w) dw + \frac{1}{2i\pi}\int_{(c')} \psi(f,s+w)^2 \Big(\sum_{n\leq Z}()\Big) Z^w \hat{\varphi}(w) dw$$

(avec $() = c_f(n) n^{s+w-1}$, $-1 < -c < -\sigma = \operatorname{Re}(s)$ et $-1 < c' < 0$.)

Pour $s = 1/2 + it$, ($|t|\leq T$), $c = -1/2 - (\log qT)^{-1}$, $c' = -(\log qT)^{-1}$, le point essentiel est que l'on a sur $\operatorname{Re}(w) = c$
$$|\psi(s+w,f)^2 (qT)^w| \ll (qT)^{1/2} T$$

et sur $\operatorname{Re}(w) = c'$
$$|\psi(s+w,f)^2 (qT)^w| \ll 1$$

qui sont aussi bonnes que les inégalités de Bombieri dans la dépendance en $q$.

Le reste est le même en utilisant la forme de $c_f(n)$ et les diverses variantes du grand crible, en particulier (70) et (71). On rencontre un terme de la forme
$$\sum_f^h \int_{-T}^T \Big|\sum_{l,m\geq 1} c(l,m)\lambda_f(l)(lm^2)^{-1/2-it}\Big|^2 dt$$

($c(l,m) = \varepsilon(m)\tau(l)\varphi(lm^2/Z)$) qui ne relève pas directement de (71) (car $\beta = 1/2$); on peut le traiter, par exemple, en estimant directement la partie $m > qT$ de la somme en $m$, en exploitant la décroissance rapide de $\varphi$ pour voir que cette partie est négligeable – ou simplement en prenant dès le départ une fonction test $\varphi$ à support compact.

On vérifie que la pondération par $(f,f)^{-1}$, qui est évidemment cruciale, est la bonne: elle intervient dans les estimations de grand crible pour les approximations de $L(f, 1/2 + it)^2$ ci-dessus, élevées au carré avant d'être intégrées en $t$.

□

On trouve ainsi pour chacune des 3 classes la même estimation (75) qui est celle du théorème 1.1.

### 7.5 Preuve du Corollaire

Nous prouvons maintenant le Corollaire 1.1 qui résulte du précédent par une application triviale de l'inégalité de Cauchy-Schwarz:

on a la majoration

$$\left(\sum_{f \in S_2(q)^+} N_f(\alpha, T)\right)^2 \leq \sum_{f \in S_2(q)^+}^h N_f^2(\alpha, T) \times \sum_{f \in S_2(q)^+} (f,f).$$

En vertu de la borne triviale

$$N_f(\alpha, T) = O(T \log(qT)),$$

le premier terme de droite est majoré par

$$O(T^B q^{\frac{3(1-\alpha)}{2-\alpha} - 1} \log^{A+1} q).$$

Le deuxième terme est majoré par $O(q \dim J_0^n(q) \log^3 q)$ en vertu du lemme suivant ([M-M] p. 337):

**Lemme 7.6** *Pour $f \in S_2(q)^+$, on a la majoration*

$$(f,f) \ll q \log^3 q.$$

D'après [Br] Lemma 3.9, on a l'encadrement

$$c_1 \phi(q) \leq \dim J_0^n(q) \leq c_2 \phi(q);$$

d'autre part $q/\phi(q) = O(\log \log q)$, on en déduit, que pour $c > 0$ assez petit et $A'$ assez grand, on a

$$\sum_{f \in S_2(q)^+} N_f(\alpha, T) = O(\dim J_0^n(q) T^B q^{-c(\alpha - 1/2)} \log q^{A'})$$

Ce qui conclut la preuve du Corollaire.

## 8 Le second théorème de densité

La preuve du théorème 1.2 est une généralisation à $GL(2)$ de la méthode des pseudo-caractères inventée par Selberg [M-S].

Comme précédemment, on cherche à construire un polynôme détecteur des zéros de $L(f, s)$. Or le fait que $L(f, s)$ est un produit eulérien de polynômes de degré 2 rend le jeu combinatoire très compliqué aussi nous avons recours à une astuce. En effet on va factoriser $L(f, s)$ par un produit eulérien de polynômes de degré 1: soit $P$ l'entier formé du produit de tous les nombres premiers $p$ tels que $1997^2/p > 1$, on a alors la factorisation suivante dans le demi-plan $\text{Re}(s) > 1$

$$(76) \qquad L(f, s) = \prod_{p|P}(1 - \lambda_f(p)p^{-s} + \varepsilon(p)p^{-2s})^{-1} G(f, s)\tilde{L}_P(f, s)$$

avec

$$\tilde{L}_P(f, s) = \sum_{\substack{n \geq 1 \\ (n, P)=1}} \mu^2(n)\lambda_f(n)n^{-s} = \sum_{n \geq 1} \tilde{\lambda}_f(n)n^{-s},$$

disons. Par la factorisation ($p \nmid q$)

$$(1 - \lambda_f(p)p^{-s} + p^{-2s})^{-1} = (1 - \alpha_p p^{-s})^{-1}(1 - \overline{\alpha_p}p^{-s})^{-1}$$

avec $|\alpha_p| = 1$, et le choix de l'entier $P$, on voit que la fonction

$$\prod_{p|P}(1 - \lambda_f(p)p^{-s} + \varepsilon(p)p^{-2s})^{-1} G(f, s)$$

est holomorphe dans le demi-plan $\text{Re}(s) > 1/2$, et est minorée et majorée uniformément dans tout domaine de la forme $\text{Re}(s) \geq 1/2 + \epsilon$, avec $\epsilon > 0$, par des constantes positives ne dépendant que de $P$ et de $\epsilon$ et donc indépendantes de $f$; c'est la raison d'être de l'entier $P$: en criblant par les facteurs premiers de $P$, on va uniformiser les comportements des différents produit eulériens qui apparaîtront dans la suite...

En particulier, les zéros de $L(f, s)$ de partie réelle $> 1/2$ sont ceux de $\tilde{L}_P(f, s)$. Comme on n'étudie que les zéros de $L(f, s)$ de partie réelle éloignée de $1/2$, on ne perd pas beaucoup à faire cette factorisation; elle va cependant simplifier considérablement l'aspect combinatoire de notre problème.

Dans la première partie nous donnons une construction des pseudo-caractères que nous utiliserons dans la suite.

### 8.1 Pseudo-caractères de formes modulaires

Soit $0 < \delta < 1/2$ un paramètre. Pour tout $f \in S_2(q)^+$, on définit un ensemble d'entiers $R_{f,\delta}$:

$$R_{f,\delta} := \{r,\ \mu^2(r) = 1,\ (r, qP) = 1,\ p|r \Rightarrow |\lambda_f(p)| > p^{-\delta}\}$$

(on abrègera $R_{f,\delta}$ en $R_f$ dans la suite).

Sur l'ensemble $R_f$, on définit la fonction arithmétique $c_f(r)$ par la formule

$$c_f(r) := \mu(r)\frac{r}{\lambda_f(r)^2}.$$

Pour chaque $r \in R_f$ on définit ensuite sur tout les entiers le pseudo-caractère

$$c_{f,r}(n) = \mu^2(n) c_f((n,r)).$$

Quand la forme $f$ sera fixée ou que le contexte le permettra, on oubliera l'indice $f$ dans l'écriture de ces fonctions. Le lemme suivant est la transposition directe du Lemma 2 de [Ju]:

**Lemme 8.1** *Soient $f, g \in S_2(q)^+$, $r \in R_f$, $r' \in R_g$; On définit les nombres $h(d; f, r; g, r')$ par l'équation*

$$\sum_{d=1}^{\infty} h(d; f, r; g, r') d^{-s}$$
$$= \prod_{\substack{p|r \\ p \nmid (r,r')}} \left(1 + (c_f(p) - 1)p^{-s}\right) \prod_{\substack{p|r' \\ p \nmid (r,r')}} \left(1 + (c_g(p) - 1)p^{-s}\right) \prod_{p|(r,r')} \left(1 + (c_f(p)c_g(p) - 1)p^{-s}\right)$$

*(en particulier la somme précédente est finie et parcourt l'ensemble des diviseurs de $[r, r']$).*
*Alors on a l'égalité*

$$c_{f,r}(n) c_{g,r'}(n) = \mu^2(n) \sum_{d|n} h(d; f, r; g, r').$$

Le lemme suivant est l'adaptation du Lemma 3. de *loc. cit.*

**Lemme 8.2** *On a l'égalité*

$$\sum_{d=1}^{\infty} \frac{h(d; f, r; f, r')}{d} \lambda_f(d)^2 \prod_{p|d} (1 + \frac{\lambda_f(p)^2}{p})^{-1} = \delta_{r,r'} |c_f(r)|,$$

*où $\delta_{r,r'}$ est le symbole de Kronecker.*

*Preuve.* — La première partie du lemme résulte de ce que la somme étudiée vaut

$$\prod_{\substack{p|rr' \\ p \nmid (r,r')}} \left(1 + \frac{(c_f(p) - 1)}{p} \lambda_f(p)^2 (1 + \frac{\lambda_f(p)^2}{p})^{-1}\right) \prod_{p|(r,r')} \left(1 + \frac{(c_f(p)^2 - 1)}{p} \lambda_f(p)^2 (1 + \frac{\lambda_f(p)^2}{p})^{-1}\right),$$

et la preuve résulte de la définition de $c_f(p)$.

□

## 8.2 Théorème de grand crible avec pseudo-caractères

Le but de cette section est de montrer le théorème suivant (une forme faible d'un résultat de grand crible avec pseudo-caractères pour les formes modulaires.)
On notera $|X_0(q)| := Vol(X_0(q)) = \frac{\pi}{3} q \prod_{p|q}(1 + p^{-1})$,

**Théorème 8.3** *Soit $1 \leq N$ et $R^2 \leq N$. Soit $(\beta_n)_{1 \leq n \leq N}$ un ensemble de complexes, alors on a la majoration suivante: pour tout $\varepsilon > 0$,*

$$\sum_{f \in S_2(q)^+}^h \sum_{\substack{r \in R_f \\ r \leq R}} \frac{1}{|c_f(r)|} \Big| \sum_{N \leq n \leq 2N} \beta_n c_{f,r}(n) \tilde{\lambda}_f(n) \Big|^2$$

$$\ll_\varepsilon \big( \frac{M}{|X_0(q)|} + q^{1/2+\varepsilon} N^{1/2+\varepsilon} R^{1+\delta+\varepsilon} \big) \sum_{N \leq n \leq 2N} |\beta_n|^2.$$

On en déduit alors le corollaire suivant à la manière de [Mon], Chap. 7.

**Corollaire 8.4** *On reprend les notations du Théorème précédent. Soit $T \geq 1$ et $1/2 < \alpha \leq 1$; pour chaque $f \in S_2(q)^+$, soit $J_f \subset R(\alpha, T)$ un ensemble fini de points bien-espacés au sens suivant:*

*pour tous $\rho, \rho' \in J_f$ distincts on a $\mathrm{Im}(\rho - \rho') \geq 1/\log q$.*

*On suppose en plus que $\beta_n = 0$ si $n \leq eR^2$. Pour tout $\varepsilon > 0$,*

$$\sum_{f \in S_2(q)^+}^h \sum_{\rho \in J_f} \sum_{\substack{r \in R_f \\ r \leq R}} \frac{1}{|c_f(r)|} \Big| \sum_{n \leq N} \beta_n c_{f,r}(n) \tilde{\lambda}_f(n) n^{-\rho} \Big|^2$$

$$\ll_\epsilon T^B (\log q + \log N)(1 + \log(\frac{\log 2N}{\log 2R})) \sum_{n \leq N} \frac{|\beta_n|^2}{n^{2\alpha}} \big( \frac{n}{|X_0(q)|} + n^{1/2+\varepsilon} q^{1/2+\varepsilon} R^{1+\delta+\varepsilon} \big).$$

### 8.3 Preuve du Théorème

Par dualité il suffit de montrer que pour toute famille $(\beta_{f,r})_{f \in S_2(q)^+, r \in R_f}$ on a la majoration

(77) $$\sum_{leqn \sim N} \Big| \sum_{f \in S_2(q)^+} \frac{1}{(f,f)^{1/2}} \sum_{r \leq R} \frac{1}{|c_f(r)|^{1/2}} \beta_{f,r} c_{f,r}(n) \tilde{\lambda}_f(n) \Big|^2$$

$$\ll \big( \frac{N}{|X_0(q)|} + N^{1/2+\varepsilon} q^{1/2+\varepsilon} R^{1+\varepsilon} \big) \sum_{f,r} |\beta_{n,r}|^2$$

On insère alors une fonction test $\phi(n/N)$ où $\phi$ est $C^\infty$ à support compact dans $[1/2, 3]$ et qui vaut 1 sur $[1, 2]$, on a donc

$$\sum_{n \sim N} |\sum_f \ldots|^2 \leq \sum_n \phi(\frac{n}{N}) |\sum_f \ldots|^2.$$

Ouvrant le carré, on est confronté à l'évaluation des sommes suivantes

$$\frac{1}{((f,f)|c_f(r)|(g,g)|c_g(r')|)^{1/2}} \sum_n \phi(\frac{n}{N}) c_{f,r}(n) c_{r',g}(n) \tilde{\lambda}_f(n) \lambda_g(n) := S(f, r; g, r')(N)$$

(les coefficients de Fourier de $g$ sont réels).

On applique le Lemme 8.1 pour obtenir

$$S(f,r;g,r')(N) = \frac{1}{\Delta_{f,r;g,r'}} \sum_{d \leq N} h(d;f,r;g,r')\lambda_f(d)\lambda_g(d) \sum_{\substack{n \\ (n,d)=1}} \mu^2(n)\tilde{\lambda}_f(n)\tilde{\lambda}_g(n)\phi(\frac{nd}{N});$$

avec $r \in R_f$, $r' \in R_g$

$$\Delta_{f,r;g,r'} = ((f,f)|c_f(r)|(g,g)|c_g(r')|)^{1/2}.$$

Notons

$$S_d(f;g)(N) := \sum_{\substack{n \\ (n,d)=1}} \mu^2(n)\tilde{\lambda}_f(n)\tilde{\lambda}_g(n)\phi(\frac{nd}{N}).$$

On considère alors la série de Dirichlet suivante (par hypothèse $(d,P) = 1$)

$$\begin{aligned}\tilde{L}_d(f \otimes g, s) &:= \sum_{\substack{n \\ (n,d)=1}} \mu^2(n)\frac{\tilde{\lambda}_f(n)\tilde{\lambda}_g(n)}{n^s} \\ &= \prod_{p|d}(1 + \frac{\lambda_f(p)\lambda_g(p)}{p^s})^{-1} \times \prod_{(p,P)=1}(1 + \frac{\lambda_f(p)\lambda_g(p)}{p^s}) \\ &:= H_d(s) \times \tilde{L}_P(f \otimes g, s)\end{aligned}$$

Comme $(d,P) = 1$, $H_d(s)$ est une fonction holomorphe dans tout demi-plan $\text{Re}(s) \geq \epsilon$ avec $\epsilon > 0$ et est majorée par $O_\epsilon(\tau(d))$ *indépendamment* de $f$ et $g$ (par le choix de l'entier $P$). La fonction $\tilde{L}_P(f \otimes g)$ admet la factorisation suivante:

$$\tilde{L}_P(f \otimes g, s) = L_P(f \otimes g, s)G_P(f \otimes g, s).$$

Dans la formule ci-dessus, $L_P(f \otimes g, s)$ est la convolution de Rankin-Selberg de $f$ et $g$ dont on a enlevé les facteurs eulériens correspondant aux nombres premiers $p$ divisant $P$:

$$L_P = L \times \prod_{p|P}(1 - \frac{1}{p^{2s}})(1 - \frac{\alpha_p\overline{\beta_p}}{p^s})(1 - \frac{\alpha_p\overline{\beta_p}}{p^s})(1 - \frac{\overline{\alpha_p}\beta_p}{p^s})(1 - \frac{\overline{\alpha_p}\beta_p}{p^s}),$$

où on a posé $\lambda_f(p) = \alpha_p + \overline{\alpha_p}$, $\lambda_g(p) = \beta_p + \overline{\beta_p}$, avec $|\alpha_p| = |\beta_p| = 1$. Et — par notre choix de l'entier $P$ — $G_P(f \otimes g)$ est une fonction holomorphe dans tout demi-plan $\text{Re}(s) \geq 1/2 + \epsilon$ avec $\epsilon > 0$, qui est majorée dans ce domaine par une constante ne dépendant que de $P$ et de $\delta$.

Par transformée de Mellin, on a

$$S_d(f;g)(N) = \frac{1}{2i\pi} \int_{(2)} \tilde{L}_d(f \otimes g, w)\hat{\phi}(w)(\frac{N}{d})^w dw.$$

Déplaçant le contour en $\text{Re}w = 1/2+\varepsilon$, et on rencontre éventuellement un pôle simple en 1 si et seulement si $f = g$ (qui correspond au pôle de $L(f \otimes f, s)$ en 1 de résidu $(f,f)/|X_0(q)|$);

en utilisant la majoration suivante valable uniformément pour $1/2 \leq \Re s \leq 3/4$ et $f, g$ primitives:
$$L(f \otimes g, s) \ll (q^2)^{\frac{1-\sigma}{2}},$$
on obtient alors l'égalité

$$\begin{aligned}
S_d(f;g)(N) &= \delta_{f,g} Res_{s=1} \tilde{L}_d(f \otimes g, s)\frac{N}{d} + O((\frac{N}{d})^{1/2+\varepsilon} q^{1/2+\varepsilon}) \\
&= \delta_{f,g} \gamma_P(f \otimes f) \prod_{p|d}(1 + \frac{\lambda_f(p)^2}{p})^{-1} \frac{(f,f)N}{|X_0(q)|d} + O_P((\frac{N}{d})^{1/2+\varepsilon} q^{1/2+\varepsilon})
\end{aligned}$$

où $\gamma_P(f \otimes f)$ est une constante bornée en fonction de $P$ mais indépendamment de $f$. On en déduit en utilisant la minoration d'Hoffstein-Lockhart $(f,f) \gg_\varepsilon q^{1-\varepsilon}$

$$S(f, r; g, r')(N) = TP + TE$$

avec

$$\begin{aligned}
TE &= O_P(\frac{N^{1/2+\epsilon} q^{-1/2+\varepsilon}}{|c_f(r)c_g(r')|^{1/2}} \prod_{p|r}(\frac{2|c_f(p)\lambda_f(p)|}{p^{1/2}} + 1) \prod_{p|r'}(\frac{2|c_g(p)\lambda_g(p)|}{p^{1/2}} + 1)) \\
&= O(R^{\delta+\varepsilon} N^{1/2+\varepsilon} q^{-1/2+\varepsilon}).
\end{aligned}$$

D'autre part on voit que
$$TP = 0 \ si \ f \neq g$$
et si $f = g$ (on utilise le fait que $N \geq R^2$)

$$TP = \frac{\gamma_P(f \otimes f)M}{|c_f(r)c_f(r')|^{1/2}|X_0(q)|} \sum_{d \leq N} \frac{h(d; f, r; f, r')\lambda_f(d)^2}{d} \prod_{p|d}(1 + \frac{\lambda_f(p)^2}{p^s})^{-1}$$

(78) $$= \frac{\gamma_P(f \otimes f)\delta_{r,r'}N}{|X_0(q)|}$$

On en déduit que la quantité (77) est majorée par

$$\ll_{P,\epsilon} \sum_{\substack{f \in S_2(q)^+ \\ r \leq R \\ r \in R_f}} \sum_{\substack{g \in S_2(q)^+ \\ r' \leq R \\ r' \in R_g}} |\beta_{f,r}\overline{\beta_{g,r'}}|(\delta_{r,r'}\delta_{f,g}\frac{N}{|X_0(q)|} + N^{1/2+\varepsilon} q^{-1/2+\varepsilon} R^{\delta+\varepsilon})$$

$$\ll_\epsilon (\frac{N}{|X_0(q)|} + N^{1/2+\varepsilon} q^{1/2+\varepsilon} R^{1+\delta+\varepsilon}) \sum_{f,r} |\beta_{f,r}|^2$$

et ceci conclut la preuve.

### 8.4 Polynômes détecteurs de zéros

On va maintenant construire un détecteur des zéros de $L(f,s)$, où plutôt de $\tilde{L}(f,s)$.

On introduit d'abord la fonction suivante qui minimise le crible à carré de Selberg. Soit $1 \leq z_1 < z_2$ on définit la fonction $\lambda_d = \mu(d) g_{z_1,z_2}(d)$ avec

$$g_{z_1,z_2}(d) = \begin{cases} 1 & si\ d \leq z_1 \\ \log(d/z_2)/\log(z_1/z_2) & si\ z_1 \leq d \leq z_2 \\ 0 & si\ d > z_2 \end{cases}$$

Le choix de ce "mollifier" est fondamental car il permettra de sauver un premier facteur $\log q$ (un deuxième sera sauvé grâce aux pseudo-caractères). En effet on a le résultat suivant du à Graham ([Gr] Lemma 9):

**Lemme 8.5** *Si $1/2 < \alpha < 1$, on a*

$$\sum_{n \leq x} (\sum_{d|n} \lambda_d)^2 n^{1-2\alpha} \ll \frac{\log(x/z_1)}{\log(z_2/z_1)} x^{2-2\alpha}.$$

Le lemme suivant construit un polynôme détecteur de zéros (Cf. Lemma 1 de [Ju]):

**Lemme 8.6** *Soit $r \in R_f$, on a l'égalité*

(79) $$\sum_{\substack{n \\ (n,P)=1}} \mu^2(n)(\sum_{d|n} \lambda_d) c_{f,r}(n) \frac{\lambda_f(n)}{n^s} = \tilde{L}_P(f,s) M_r(f,s);$$

*où on a posé*

$$M_r(f,s) = \sum_{\substack{d=1 \\ (d,P)=1}}^{\infty} \lambda_d \lambda_f(d) \frac{c_{f,r}(d)}{d^s} \prod_{p|r/(r,d)} (1 + \frac{\lambda_f(p) c_f(p)}{p^s}) \prod_{p|dr} (1 + \frac{\lambda_f(p)}{p^s})^{-1}.$$

*Preuve.* — Utilisant la multiplicativité de $c_{f,r}(-)\lambda_f(-)$, on voit que le terme de gauche de (79) vaut

$$= \sum_{\substack{d \\ (d,P)=1}} \lambda_d c_{f,r}(d) \frac{\lambda_f(d)}{d^s} \sum_{\substack{n \\ (n,dP)=1}} \mu^2(n) c_{f,r}(n) \frac{\lambda_f(n)}{n^s}$$

$$= \sum_{\substack{d \\ (d,P)=1}} \lambda_d c_{f,r}(d) \frac{\lambda_f(d)}{d^s} \prod_{\substack{p \\ (p,dP)=1}} (1 + c_{f,r}(p) \frac{\lambda_f(p)}{p^s})$$

$$= \prod_{\substack{p \\ (p,P)=1}} (1 + c_{f,r}(p) \frac{\lambda_f(p)}{p^s}) \sum_{\substack{d \\ (d,P)=1}} \lambda_d c_{f,r}(d) \frac{\lambda_f(d)}{d^s} \prod_{p|d} (1 + c_{f,r}(p) \frac{\lambda_f(p)}{p^s})^{-1}$$

$$= \tilde{L}(f,ss) \prod_{p|r} \Big(\frac{1 + c_f(p)\lambda_f(p)p^{-s}}{1 + \lambda_f(p)p^{-s}}\Big) \sum_{\substack{d \\ (d,P)=1}} \lambda_d c_{f,r}(d) \frac{\lambda_f(d)}{d^s} \times$$

$$\prod_{p|(r,d)}(1+c_f(p)\frac{\lambda_f(p)}{p^s})^{-1}\prod_{p|d/(r,d)}(1+\frac{\lambda_f(p)}{p^s})^{-1}$$

$$=\tilde{L}_P(f,s)M_r(f,s).$$

□

Le lemme suivant détecte les zéros de $L(f,s)$:

**Lemme 8.7** *Soit $r \in R_f$ avec $r \leq R$. Pour tout $1-1/13 \leq \alpha \leq 1$ et $T \geq 1$, soit $\rho \in R(\alpha, T)$ un zéro de $L(f,s)$. On suppose que l'on a les inégalités suivantes: pour un certain $\epsilon > 0$,*

(80) $\quad X \geq (q^{1/2}TR^{1+2\delta}z_2)^{1/(2\alpha-1)+\epsilon}, \; \log X \ll \log qT, \; \log^{1/2} q \leq \log R \leq \frac{1}{2}\log x$

*où on a posé $x = X\log^2(qT)$ et soit*

$$g_r(s,f) := \sum_{\substack{1<n\leq x \\ (n,P)=1}} \mu^2(n)(\sum_{d|n}\lambda_d)c_{f,r}(n)e^{-n/X}\frac{\lambda_f(n)}{n^s}.$$

*Alors, on a la minoration*

$$1 \ll_{P,\epsilon} |g_r(\rho;f)|;$$

*la constante impliquée ne dépendant que de $P$, de $\epsilon$ et de $\alpha'$.*

*Preuve.* — Soit $0 < \epsilon$ assez petit, par déformation d'un contour et le lemme précédent on a

$$1+\sum_{\substack{1<n \\ (n,P)=1}}\mu^2(n)(\sum_{d|n}\lambda_d)c_{f,r}(n)e^{-n/X}\frac{\lambda_f(n)}{n^\rho} = \int_{(-\alpha+1/2-\epsilon)}\tilde{L}(f,s+\beta)\Gamma(s)X^s M_r(f,s+\beta)ds.$$

Par le principe de Phragmen-Lindelöf, on a pour $0 \leq \sigma \leq 1$, la majoration

$$L(f,\sigma+it) \ll_\epsilon (q^{1/2}(|t|+1))^{1-\sigma+\epsilon}.$$

La factorisation (76), implique que cette majoration est aussi valable pour $\tilde{L}_P(f,s)$ quand $\mathrm{Re}(s) \geq 1/2+\epsilon$. D'autre part on a,

$$|M_r(f,1/2+\beta-\alpha+\epsilon+it)| \ll_\epsilon r^{1/2+2\delta+\epsilon}z_2^{1/2+2\epsilon}.$$

En effet, comme $r \in R_f$, on a

$$|\lambda_f(d)c_{r,f}(d)| \leq \tau(d)(d,r)^{1+\delta}, \quad \prod_{p|r/(r,d)}|1+\frac{\lambda_f(p)c_f(p)}{p^{1/2}}| \ll_\epsilon r^{1/2+\delta+\epsilon};$$

enfin comme $(dr,P)=1$ on a

$$\prod_{p|dr}(1+\frac{\lambda_f(p)}{p^{1/2+\epsilon}})^{-1} \leq \tau(dr).$$

On en déduit que l'intégrale est $O_\epsilon((q^{1/2}TR^{1+2\delta}z_2)^{1/2+\epsilon}/X^{\alpha-1/2})$. Enfin, on voit facilement que si la série est tronquée à $x = X\log^2 qT$ le reste est un $o_\epsilon(1)$.

$\square$

### 8.5 Comptage des zéros

On reprend la méthode de comptage classique utilisée par exemple dans [Ju].

On subdivise le rectangle $R(\alpha, T)$ en rectangles de la formes $R_k := [\alpha, 1] \times [k/\log q, (k+1)/\log q]$. Pour chaque $f \in S_2(q)^+$, on choisit, pour chaque $k$ un zéro $\rho_{f,k}$ de $L(s,f)$ contenu dans $R_k$; enfin réunissant les $k$ pairs et les $k$ impairs, on obtient deux ensembles de zéros $J_f^+$ et $J_f^-$ chacuns $1/\log q$-bien espacés. On commence par majorer les quantités

$$|J^+| := \sum_{f \in S_2(q)^+} |J_f^+|, \ |J^-| := \sum_{f \in S_2(q)^+} |J_f^-|.$$

Par le lemme 8.7 et le Corollaire 8.4 on a, pour peu que les conditions (80) soient remplies et que $z_1 \geq R^2$, la majoration

$$\sum_{f \in S_2(q)^+} \sum_{r \leq R} \frac{1}{|c_f(r)|} \frac{|J_f^\pm|}{(f,f)} \ll \sum_{f,r} \frac{1}{|c_f(r)|} \sum_{\rho \in J_f^\pm} \frac{|g_f(\rho; r)|^2}{(f,f)}$$

(81)
$$\ll T^B \log q \sum_{z_1 < n \leq x} (\sum_{d|n} \lambda_d)^2 (\frac{n^{1-2\alpha}}{|X_0(q)|} + n^{1/2-2\alpha}q^{3/2}R^{1+\varepsilon})$$

$$\ll T^B \log q(\frac{x^{2(1-\alpha)}}{|X_0(q)|} + z_1^{-1/2+\varepsilon}q^{3/2}R^{1+\varepsilon}).$$

Dans la dernière étape, on a utilisé le lemme 8.5, qui permet de sauver un facteur $\log x$ sous l'hypothèse $\log x \ll \log(z_2/z_1)$.

D'autre part part les méthodes de factorisation que ont été décrites (remarquer que le produit $\prod_p(1-p^{1+\delta})$ est absolument convergent), on a la minoration suivante que permet de sauver un deuxième facteur $\log q$:

$$\sum_{\substack{r \leq R \\ r \in R_f}} \frac{1}{|c_f(r)|} = \sum_{\substack{r \leq R \\ r \in R_f}} \frac{\lambda_f(r)^2}{r} \geq \sum_{\substack{r \in R_f \\ R^{1-\varepsilon} \leq r \leq R}} \frac{\lambda_f(r)^2}{r}$$

$$\gg_\delta \varepsilon \frac{(f,f)}{|X_0(q)|} \log(R^\varepsilon) + O_\varepsilon(R^{-1/2+2\varepsilon}q^{1/2+\varepsilon})$$

(82)
$$\gg_\delta \frac{(f,f)}{|X_0(q)|} \log R$$

(83)

si $R \gg_\delta q^{1+4\epsilon}$.

On obtient donc réunissant (81) et (82) la majoration

$$|J^\pm| \ll T^B(\frac{\log q}{\log R}x^{2(1-\alpha)} + |X_0(q)|q^{1/2+\varepsilon}R^{1+\delta+\varepsilon}z_1^{-1/2}).$$

On choisit alors $\delta = \epsilon$,

$$R = q^{1+4\varepsilon},\ z_2 = z_1^{1+\varepsilon},\ z_1 = q^{5+100\varepsilon}, x = (q^{13/2}T)^{1/(2\alpha-1)+200\varepsilon},$$

pour obtenir pour tout $\varepsilon > 0$ la majoration

$$|J^\pm| \ll_\varepsilon T^B(q^{13+\varepsilon})^{\frac{1-\alpha}{2\alpha-1}}).$$

On conclut alors la preuve du théorème 1.2 grâce au lemme de densité suivant

**Lemme 8.8** *Le nombre de zéros $\rho = \beta + i\tau$ de $L(f,s)$ contenus dans le carré*

$$\alpha \leq \sigma \leq 1,\ |t - \tau| \leq \frac{1}{2}(1-\alpha)$$

$$\text{est } \ll (1-\alpha)\log(q(|t|+1) + 1.$$

En effet, le facteur $(1-\alpha)\log(qT)$ est absorbé par le facteur $T^B q^{\varepsilon\frac{1-\alpha}{2\alpha-1}}$.

# Références

E. KOWALSKI: ekowalsk@math.rutgers.edu
Dept. of Math. RUTGERS University, New Brunswick, NJ 08903 USA

P. MICHEL: michel@math.u-psud.fr
Université PARIS-SUD, Bât. 425 91405 ORSAY Cdx. FRANCE